\documentclass[aos,MSNbibl,seceqn,nameyear,dvips]{arximspdf}
\usepackage{dcolumn}
\usepackage{graphicx}

\doi{10.1214/11-AOS905}
\volume{39}
\issue{5}
\pubyear{2011}
\firstpage{2410}
\lastpage{2447}

\makeatletter

\newcolumntype{d}[1]{D{.}{.}{#1}}

\newtheorem{theorem}{Theorem}
\newtheorem{prop}{Proposition}

\newproclaim{rem}{Remark}

\makeatother

\begin{document}
\begin{frontmatter}

\title{Factor models and variable selection in high-dimensional regression
analysis}
\runtitle{Factor models and variable selection}

\begin{aug}
\author[A]{\fnms{Alois} \snm{Kneip}\corref{}\ead[label=e1]{akneip@uni-bonn.de}}
\and
\author[B]{\fnms{Pascal} \snm{Sarda}\ead[label=e2]{Pascal.Sarda@math.ups-tlse.fr}}
\runauthor{A. Kneip and P. Sarda}
\affiliation{Universit\"at Bonn and Universit\'e Paul Sabatier}
\address[A]{Statistische Abteilung\\
Department of Economics\\
\quad and Hausdorff\\
Center for Mathematics\\
Universit\"{a}t Bonn\\
Adenauerallee 24-26\\
53113 Bonn\\
Germany\\
\printead{e1}}
\address[B]{Institut de Math\'ematiques\\
Laboratoire de Statistique\\
\quad et Probabilit\'es\\
Universit\'e Paul Sabatier\\
UMR 5219 \\
118, Route de Narbonne\\
31062 Toulouse Cedex\\
France\\
\printead{e2}} 
\end{aug}

\received{\smonth{3} \syear{2011}}

%
\begin{abstract}
The paper considers linear regression problems where the number of
predictor variables is possibly larger than the sample size. The
basic motivation of the study is to combine the points of view of model
selection and functional regression by using a factor approach: it is
assumed that the predictor vector can be decomposed into a sum of two
uncorrelated random components reflecting common factors and specific
variabilities of the explanatory variables. It is shown that the
traditional assumption of a sparse vector of parameters is restrictive
in this context. Common factors may possess a significant influence on
the response variable which cannot be captured by the specific effects
of a small number of individual variables. We therefore propose to
include principal components as additional explanatory variables in an
augmented regression model. We give finite sample inequalities for
estimates of these components. It is then shown that model selection
procedures can be used to estimate the parameters of the augmented
model, and we derive theoretical properties of the estimators. Finite
sample performance is illustrated by a simulation study.
\end{abstract}

%

\end{frontmatter}

\section{Introduction}\label{sec1}

The starting point of our analysis is a high-dimensional linear
regression model
of the form
%
%
\begin{equation}\label{linear-model}
Y_i = \bolds{\beta}^T\mathbf{X}_i + \varepsilon_i,\qquad i=1,\ldots,n,
\end{equation}
where $(Y_i,\mathbf{X}_i)$, $i=1,\ldots,n$, are i.i.d. random pairs
with $Y_i\in\mathbb{R}$ and $\mathbf{X}_i=(X_{i1},\ldots,X_{ip})^T\in
\mathbb{R}^p$. We will assume without loss of generality that
$\mathbb{E}(X_{ij})=0$ for all $j=1,\ldots,p$.
Furthermore,\vspace*{1pt} $\bolds{\beta}$ is a vector of parameters in
$\mathbb{R}^p$ and $(\varepsilon_i)_{i=1,\ldots,n}$ are centered i.i.d.
real random variables independent with $\mathbf{X}_i$ with
$\operatorname{Var}(\varepsilon_i)=\sigma^2$. The dimension $p$ of the
vector of parameters is assumed to be typically larger than the sample
size $n$.

Roughly speaking, model (\ref{linear-model}) comprises two main
situations which have been considered independently in two separate
branches of statistical literature. On one side, there is the situation
where $\mathbf{X}_i$ represents a (high-dimensional) vector of
different predictor variables. Another situation arises when the
regressors are $p$ discretizations (e.g., at different observations
times) of a same curve. In this case
model (\ref{linear-model}) represents a discrete version of an
underlying continuous \textit{functional linear model}. In the two
setups, very different strategies for estimating
$\bolds{\beta}$ have been adopted, and underlying structural
assumptions seem to
be largely incompatible. In this paper we will study similarities and
differences of
these methodologies, and we will show that
a combination of ideas developed in the two settings leads to new
estimation procedures
which may be useful in a number of important applications.

The first situation is studied in a large literature on model selection
in high-dimensional
regression. The basic structural assumptions can be described as follows:
\begin{itemize}
\item There is only a relatively small number of predictor variables
with $|\beta_j|>0$ which have a significant
influence on the outcome $Y$. In other words,
the set of nonzero coefficients is sparse, $S:=\#\{j|\beta_j\neq0\}\ll p$.
\item The correlations between different explanatory variables $X_{ij}$
and $X_{il}$, $j\neq l$, are ``sufficiently'' weak.
\end{itemize}
The most popular procedures to identify and estimate nonzero
coefficients $\beta_j$ are Lasso and the Dantzig selector. Some
important references are \citet{Ti96}, \citet{MeBu06},
\citet{ZhYu}, \citet{Va08}, \citet{BiRiTss09},
\citet{CaTa07} and \citet{Ko09}. Much work in this domain is
based on the assumption that the columns $(X_{1j},\ldots, X_{nj})^T$,
$j=1,\ldots,p$, of the design matrix are almost orthogonal. For example,
\citet{CaTa07} require that ``every set of columns with
cardinality less than $S$ approximately behaves like an orthonormal
system.'' More general conditions have been introduced by
\citet{BiRiTss09} or \citet{ZhVaBu09}. The theoretical
framework developed in
these papers also allows one to study model selection for regressors
with substantial amount of correlation, and it provides a basis for the
approach presented in our paper.

In sharp contrast, the setup considered in the literature on functional
regression rests upon
a very different type of structural assumptions. We will consider the
simplest case
that $X_{ij}=X_i(t_j)$ for random functions $X_i\in L^2([0,1])$
observed at an equidistant
grid $t_j=\frac{j}{p}$. Structural assumptions on coefficients and
correlations between
variables can then be subsumed as follows:
\begin{itemize}
\item$\beta_j:=\frac{\beta(t_j)}{p}$, where $\beta(t)\in L^2([0,1])$
is a continuous slope function,
and as $p\rightarrow\infty$, $\sum_j \beta_j X_{ij}=\sum_j \frac{\beta
(t_j)}{p} X_i(t_j)
\rightarrow\int_0^1 \beta(t)X_i(t)\,dt$.
\item There are very high correlations between
explanatory variables $X_{ij}=X_i(t_j)$ and $X_{il}=X_i(t_l)$, $j\neq
l$. As $p\rightarrow\infty$, $\operatorname{corr}(X_i(t_j),X_i(t_{j+m}))\rightarrow
1$ for any fixed~$m$.
\end{itemize}
Some important applications as well as theoretical results on
functional linear regression are, for example, presented in
\citet{RaDa91}, \citet{CaFeSa99}, \citet{CuFeFr02},
\citet{YaMuWa05}, \citet{CaHa06}, \citet{HaHo07},
\citet{CaMaSa07} and \citet{CrKnSa09}. Obviously, in this
setup no variable
$X_{ij}=X_i(t_j)$ corresponding to a specific observation at grid point
$t_j$ will possess a particulary high influence on $Y_i$, and there will
exist a large number of small, but nonzero coefficients $\beta_j$ of
size proportional to $1/p$. One may argue that dimensionality reduction
and therefore some underlying concept of ``sparseness'' is always
necessary when dealing with high-dimensional problems. However, in
functional regression sparseness is usually not assumed with respect to
the coefficients $\beta_j$, but the model is rewritten using a
``sparse'' expansion of the predictor functions~$X_i$.\looseness=-1

The basic idea relies on the so-called Karhunen--Lo\`eve decomposition
which provides a decomposition of
random functions in terms of functional principal components of the
covariance operator of $X_i$.
In the discretized case analyzed in this paper this amounts to consider
an approximation of $X_i$ by
the principal components
of the covariance matrix $\bolds{\Sigma}=\mathbb{E}(\mathbf
{X}_i\mathbf{X}_i^T)$. In practice, often a small number $k$ of
principal components will suffice
to achieve a small $L^2$-error. An important points is now that even if
$p>n$ the eigenvectors corresponding to the leading eigenvalues $\mu
_1,\ldots,\mu_k$ of $\bolds{\Sigma}$
can be well estimated by the
eigenvectors (estimated principal
components) $\widehat{\bolds{\psi}}_r$ of the empirical covariance
matrix $\widehat{\bolds{\Sigma}}$.
This is due to the fact that if the predictors
$X_{ij}$ represent discretized values of a continuous functional
variable, then for sufficiently
small $k$ the eigenvalues $\mu_1,\ldots,\mu_k$ will necessarily be of an
order larger
than $\frac{p}{\sqrt{n}}$ and will thus exceed the magnitude of purely
random components.
From a more general point of view the underlying theory will be
explained in detail in Section \ref{sec4}.

Based on this insight, the most frequently used approach in functional
regression is to approximate $\mathbf{X}_i\approx\sum_{r=1}^k
(\widehat{\bolds {\psi}}^T_r\mathbf{X}_i) \widehat{\bolds{\psi}}_r $ in
terms\vspace*{2pt} of the first $k$ estimated principal\vspace*{1pt}
components $\widehat{\bolds{\psi}}_1,\ldots, \widehat{\bolds{\psi}}_k$,
and to rely on the approximate model $Y_i\approx\sum_{r=1}^k \alpha_r
\widehat{\bolds{\psi}}{}^T_r\mathbf {X}_i+\varepsilon_i$. Here, $k$
serves as smoothing parameter. The new coefficients $\bolds{\alpha}$
are estimated by least squares, and $\widehat{\beta}_j=\sum_{r=1}^k
\widehat{\alpha}_r \widehat{\psi}_{rj}$. Resulting rates of convergence
are given in \citet{HaHo07}.

The above arguments show that a suitable regression analysis will
have to take into account the underlying structure of the explanatory
variables $X_{ij}$.
The basic motivation of this paper now is to combine the points of view of
the above branches of literature
in order to develop a new approach for model
adjustment and variable selection
in the practically
important situation of strongly correlated regressors. More precisely,
we will concentrate
on \textit{factor models} by assuming
that the $\mathbf{X}_i\in\mathbb{R}^p$ can be decomposed in the form
%
%
\begin{equation}\label{predictor}
\mathbf{X}_i = \mathbf{W}_i+\mathbf{Z}_i,\qquad i=1,\ldots,n,
\end{equation}
where $\mathbf{W}_i$ and $\mathbf{Z}_i$ are two uncorrelated random
vectors in $\mathbb{R}^p$. The random vector $\mathbf{W}_i$ is intended
to describe high correlations of the $X_{ij}$ while the components
$Z_{ij}$, $j=1,\ldots,p$, of $\mathbf{Z}_i$ are uncorrelated. This implies
that the covariance matrix $\bolds{\Sigma}$ of $\mathbf{X}_i$
adopts the decomposition
%
%
\begin{equation}\label{var-decomp}
\bolds{\Sigma} = \bolds{\Gamma} + \bolds{\Psi},
\end{equation}
where $\bolds{\Gamma}=\mathbb{E}(\mathbf{W}_i\mathbf{W}_i^T)$, while
$ \bolds{\Psi}$ is a diagonal matrix with diagonal entries
$\operatorname{var}(Z_{ij})$, $j=1,\ldots,p$.

Note that factor models can be found in any textbook on multivariate
analysis and must be
seen as one of the major tools
in order to analyze samples of high-dimensional vectors.
Also recall that a standard factor model is additionally based on the
assumption that
a finite
number $k$ of factors suffices to approximate $\mathbf{W}_i$ precisely. This
means that the matrix $\bolds{\Gamma}$ only possesses $k$ nonzero
eigenvalues.
In the following we will more generally assume that a small number of
eigenvectors of
$\bolds{\Gamma}$ suffices to approximate $\mathbf{W}_i$ with high accuracy.

We want to emphasize that the typical structural assumptions to be
found in the literature
on high-dimensional regression are special cases of
(\ref{predictor}). If $\mathbf{W}_i=0$ and thus $\mathbf{X}_i=\mathbf
{Z}_i$, we are in the
situation of uncorrelated regressors which has been widely studied in
the context of model
selection. On the other hand, $\mathbf{Z}_i=0$ and thus $\mathbf
{X}_i=\mathbf{W}_i$
reflect the structural assumption of functional regression.

In this paper we assume that $W_{ij}$ as well as $Z_{ij}$ represent
nonnegligible parts of the
variance of $X_{ij}$. We believe that this approach may well describe
the situation encountered in many relevant applications. Although
standard factor models are usually considered in the case
$p\ll n$, (\ref{predictor}) for large values of $p$ may be of
particular interest in time
series or spatial analysis. Indeed, factor models for large $p$ with a
finite number $k$
of nonzero eigenvalues of $\bolds{\Gamma}$ play an important role
in the
econometric study of multiple time series and panel data. Some
references are \citet{FoLi97}, \citet{Foal00}, \citet
{StWa02}, \citet{BeBo03} and
Bai (\citeyear{Bai03}, \citeyear{Bai09}).

Our objective now is to study linear regression (\ref{linear-model})
with respect to
explanatory variables which adopt decomposition (\ref{predictor}). Each
single variable $X_{ij}$,
$j=1,\ldots,p$, then possesses a \textit{specific} variability induced by
$Z_{ij}$ and may
thus explain some part of the outcome $Y_i$. One will, of course,
assume that only few
variables have a significant influence on $Y_i$ which enforces the use
of model selection
procedures.

On the other hand, the term $W_{ij}$ represents a \textit{common}
variability. Corresponding principal components quantify
a simultaneous variation of many individual regressors. As a
consequence, such principal
components may possess some additional power for predicting $Y_i$
which may
go beyond the effects of individual variables. A~rigorous discussion
will be given in Section \ref{sec3}. We want to note that the concept of
``latent variables,'' embracing the common influence of a large group
of individual
variables, plays a prominent role in applied, parametric multivariate analysis.

These arguments motivate the main results of this paper. We propose to
use an ``augmented''
regression model which includes principal components as additional
explanatory variables.
Established model selection procedures like the Dantzig selector or the
Lasso can then be
applied to estimate the nonzero coefficients of the augmented model. We
then derive theoretical
results providing bounds for the accuracy of the resulting estimators.

The paper is organized as follows: in Section \ref{sec2} we formalize
our setup.
We show in Section \ref{sec3} that the traditional sparseness
assumption is restrictive
and that a valid model may have to include principal components. The
augmented model is thus introduced with an estimation procedure.
Section \ref{sec4} deals with
the problem how accurately true principal components can be estimated
from the sample
$ \mathbf{X}_1,\ldots,\mathbf{X}_n$. Finite sample inequalities are
given, and we show
that it is possible to obtain sensible estimates of those components
which explain a considerable percentage of the
total variance of all $X_{ij}$, $j=1,\ldots,p$. Section \ref{sec5}
focuses on
theoretical properties of the augmented model, while in Section \ref
{sec6} we present
simulation results illustrating the
finite sample performance of our estimators.

\section{The setup}\label{sec2}

We study regression of a response
variable $Y_i$ on a set of i.i.d. predictors $\mathbf{X}_i\in\mathbb{R}^p$,
$i=1,\ldots,n$, which adopt
decomposition (\ref{predictor}) with $\mathbb{E}(X_{ij})=\mathbb{E}(W_{ij})
=\mathbb{E}(Z_{ij})=0$, $\mathbb{E}(Z_{ij}Z_{ik})=0$, $\mathbb
{E}(W_{ij}Z_{il})=0$, $\mathbb{E}(Z_{ij}Z_{ik}Z_{il}Z_{im})=0$ for
all $j,k,l,m\in\{1,\ldots,p\}$, $j\notin \{k,l,m\}$.
Throughout the following sections we additionally assume that there
exist constants
$D_0,D_3<\infty$ and $0<D_1\leq D_2<\infty$ such that with
$\sigma_j^2:=\operatorname{Var}(Z_{ij})$ the following assumption (A.1) is satisfied
for all~$p$:
\begin{longlist}[(A.1)]
\item[(A.1)]
$0 < D_1\leq\sigma_j^2 \leq D_2, \mathbb{E}(X_{ij}^2)\leq D_0,
\mathbb{E}(Z_{ij}^4)\leq D_3$
for all $j=1,\ldots,p$.
\end{longlist}
Recall that $\bolds{\Sigma}=\mathbb{E}(\mathbf{X}_i\mathbf
{X}_i^T)$ is the covariance matrix of $\mathbf{X}_i$ with $\bolds
{\Sigma}=\bolds{\Gamma}+\bolds{\Psi}$, where $\bolds
{\Gamma}=\mathbb{E}(\mathbf{W}_i\mathbf{W}_i^T)$ and $\bolds{\Psi
}$ is a diagonal matrix with diagonal entries $\sigma^2_j$, $j=1,\ldots
,p$. We denote as $\widehat{\bolds{\Sigma}}=\frac{1}{n}\sum
_{i=1}^n\mathbf{X}_i\mathbf{X}_i^T$ the empirical covariance matrix
based on the sample $\mathbf{X}_i$, $i=1,\ldots,n$.

Eigenvalues and eigenvectors of the standardized matrices $\frac
{1}{p}\bolds{\Gamma}$ and $\frac{1}{p}\bolds{\Sigma}$ will
play a central
role. We will use $\lambda_1\geq\lambda_2\geq\cdots$ and $\mu_1\geq\mu
_2\geq
\cdots$ to denote the eigenvalues of $\frac{1}{p}\bolds{\Gamma}$
and $\frac{1}{p}\bolds{\Sigma}$, respectively, while $\bolds
{\psi}_1,\bolds{\psi}_2,\ldots$ and $\bolds{\delta
}_1,\bolds{\delta}_2,\ldots$ denote corresponding
orthonormal eigenvectors.
Note that all eigenvectors of $\frac{1}{p}\bolds{\Sigma}$ and
$\bolds{\Sigma}$ (or
$\frac{1}{p}\bolds{\Gamma}$ and
$\bolds{\Gamma}$) are identical, while
eigenvalues differ by the factor $1/p$. Standardization is important
to establish convergence results
for large $p$, since the largest eigenvalues of $\Sigma$ tend to
infinity as $p\rightarrow\infty$.

From a conceptional point of view we will concentrate on the case that
$p$ is large
compared to $n$. Another crucial, \textit{qualitative} assumption
characterizing our approach is the dimensionality reduction of $\mathbf
{W}_i$ using a small number $k\ll p$ of
eigenvectors (principal components) of
$\frac{1}{p}\bolds{\Gamma}$ such that (in a good approximation)
$\mathbf{W}_i
\approx\sum_{r=1}^k \xi_{ir} \bolds{\psi}_r$.
We also assume that $\mathcal{D}_X=\frac{1}{p}\sum_{j=1}^p \mathbb{E}(X_{ij}^2)>
\mathcal{D}_W=\frac{1}{p}\sum_{j=1}^p \mathbb{E}(W_{ij}^2)
\gg\frac{1}{p}$.
\textit{Then all leading principal components of $\frac{1}{p}\bolds
{\Gamma}$
corresponding to the $k$ largest eigenvalues explain a considerable percentage
of the total variance of $\mathbf{W}_i$ and $\mathbf{X}_i$.}

Indeed, if $\mathbf{W}_i = \sum_{r=1}^k \xi_{ir} \bolds{\psi}_r$, we
necessarily have $\lambda_1\geq\frac{\mathcal{D}_W}{k}\gg\frac{1}{p}$
and $\mu_1\geq \lambda_1\geq\frac{\mathcal{D}_W}{k}\gg\frac{1}{p}$.
Then $\operatorname{tr}(\frac
{1}{p}\bolds{\Gamma})=\sum_{r=1}^p\lambda_r =\frac{1}{p}\sum_{j=1}^p
\mathbb{E}(W_{ij}^2)$, and the first principal component of
$\frac{1}{p}\bolds{\Gamma}$ explains a considerable proportion
$\frac{\lambda_1}{(1/p)\sum_{j=1}^p \mathbb{E}(W_{ij}^2)}\geq
\frac{1}{k}\gg\frac{1}{p}$ of the total variance of $\mathbf{W}_i$.

We want to emphasize that this situation is very different from the
setup which is usually considered in the literature on the analysis of
high-dimensional covariance matrices; see, for example, \citet
{BiLe08}. It is then assumed that the variables of interest are only
weakly correlated and that the largest\vspace*{1pt} eigenvalue $\mu_1$
of the corresponding scaled covariance matrix
$\frac{1}{p}\bolds{\Sigma}$ is of order $\frac {1}{p}$. This means that
for large $p$ the first principal component only\vspace*{1pt} explains
a negligible percentage of the total variance of $\mathbf{X}_i$,
$\frac{\mu_1}{(1/p)\sum _{j=1}^p \mathbb{E}(X_{ij}^2)}=O(\frac{1}{p})$.
It is well known that in this case no consistent estimates of
eigenvalues and principal components can be obtained from an
eigen-decomposition of
$\frac{1}{p}\widehat{\bolds{\Sigma}}$.\vspace*{1pt}

However, we will show in Section \ref{sec4}
that principal components which are able to explain a considerable
proportion of total variance can be estimated consistently. These
components will be an intrinsic part of the augmented model presented
in Section \ref{sec3}.

We will need a further assumption which ensures that all covariances
between the different variables are well approximated by their
empirical counterparts:
\begin{longlist}[(A.2)]
\item[(A.2)]
There exists a $C_0<\infty$ such that
%
%
\begin{eqnarray}
\label{B21}
\sup_{1\leq j, l\leq p}\Biggl| \frac{1}{n}\sum
_{i=1}^nW_{ij}W_{il}-\operatorname{cov}(W_{ij},W_{il})\Biggr| &\leq&
C_0\sqrt{ \frac{\log
p}{n}}, \\
\label{B22}
\sup_{1\leq j, l\leq p}\Biggl| \frac{1}{n}\sum_{i=1}^n
Z_{ij}Z_{il}-\operatorname{cov}(Z_{ij},Z_{il})\Biggr|
&\leq& C_0 \sqrt{\frac{\log p}{n}},\\
\label{B23}
\sup_{1\leq j, l\leq p}\Biggl| \frac{1}{n}\sum_{i=1}^n
Z_{ij}W_{il}\Biggr| &\leq& C_0 \sqrt{\frac{\log p}{n}},\\
\label{B24}
\sup_{1\leq j, l\leq p}\Biggl| \frac{1}{n}\sum
_{i=1}^nX_{ij}X_{il}-\operatorname{cov}(X_{ij},X_{il})\Biggr|
&\leq& C_0 \sqrt{\frac{\log
p}{n}}
\end{eqnarray}
\end{longlist}
hold simultaneously with probability $A(n,p)>0$, where
$A(n,p)\rightarrow1$ as $n,p\rightarrow\infty$, $\frac{\log
p}{n}\rightarrow0$.

The following proposition provides a general sufficient condition on
random vectors for which (A.2) is satisfied provided that the rate of
convergence of $\frac{\log p}{n}$ to 0 is sufficiently fast.
\begin{prop}\label{propcov}
Consider independent and identically distributed random vectors $\mathbf
{V}_i\in\mathbb{R}^p$, $i=1,\ldots,n$, such that for $j=1,\ldots,p$,
$\mathbb{E}(V_{ij})=0$ and
%
%
\begin{equation}\label{exponential-condition}
\mathbb{E}\bigl(e^{a|V_{ij}|}\bigr) \le C_1
\end{equation}
for positive constants $a$ and $C_1$ with moreover $\mathbb
{E}(V_{ij}^4)\le C_1$. Then, for any positive constant $C_0$ such that
$C_1^{1/2}\le\frac{1}{2}\sqrt{\frac{C_0n}{\log p}}$ and $C_1\le\frac
{1}{8}C_0e^{a\sqrt{{C_0n}/({\log p})}}\sqrt{\frac{\log p}{n}}$
%
%
\begin{eqnarray}\label{propcovresult}
&&
P\Biggl(\sup_{1\leq j, l\leq p}\Biggl| \frac{1}{n}\sum_{i=1}^n
V_{ij}V_{il}-\operatorname{cov}(V_{ij},V_{il})\Biggr| \leq C_0 \sqrt{\frac{\log p}{n}}
\Biggr)\nonumber\\[-8pt]\\[-8pt]
&&\qquad\ge 1-p^{2-C^2_0/(8(C_1+C_0^{3/2}/3))}+2p^2nC_1e^{-({a}/{\sqrt
{2}})(C_0n/\log p)^{1/4}}.\nonumber
\end{eqnarray}
\end{prop}

Note that as $n,p\rightarrow\infty$ the right-hand side of (\ref
{propcovresult}) converges
to 1 provided that $C_0$ is chosen sufficiently large and that
$p/e^{n^{1-\tau}}=O(1)$ for some $4/5<\tau<1$. Therefore, assumption
(A.2) is satisfied if the components of the random variables
$X_{ij}$ possess some exponential moments. For the specific case of
centered normally distributed random variables, a more precise bound in
(\ref{propcovresult}) may be obtained using Lemma 2.5 in \citet{ZhVaBu09}
and large deviations inequalities obtained by \citet{ZhLaWa08}. In
this case it may also be shown
that for sufficiently large $C_0$ events (\ref{B21})--(\ref{B24}) hold
with probability tending to 1 as $p\rightarrow\infty$ without any
restriction on the quotient $\log p/n$.
Of course, the rate $\sqrt{\frac{\log p}{n}}$ in (\ref{B21})--(\ref
{B24}) depends on the tails of the distributions: it would be possible
to replace this rate with a slower one in case of heavier tails than in
Proposition \ref{propcov}. Our theoretical results could be modified
accordingly.

\section{The augmented model}\label{sec3}

Let us now consider the structural model (\ref{predictor}) more
closely. It implies that
the vector $\mathbf{X}_i$ of predictors can be decomposed into
two uncorrelated random vectors $\mathbf{W}_i$ and $\mathbf{Z}_i$.
Each of these
two components separately may possess a significant influence on the
response variable $Y_i$. Indeed, if $\mathbf{W}_i$ and $\mathbf{Z}_i$
were known, a possibly substantial improvement of model (\ref
{linear-model}) would consist
in a regression of $Y_i$ on the $2p$ variables $\mathbf{W}_i$ \textit{and}~$\mathbf{Z}_i$
%
%
\begin{equation}\label{eq-30}
Y_i =
\sum_{j=1}^p\beta_j^*W_{ij}+\sum_{j=1}^p\beta_jZ_{ij}+\varepsilon_i,\qquad
i=1,\ldots,n,
\end{equation}
with different sets of parameters $\beta_j^*$ and $\beta_j$, $j=1,\ldots
,p$, for each contributor. We here again assume that $\varepsilon_i$,
$i=1,\ldots,n$, are
centered i.i.d. random variables with $\operatorname{Var}(\varepsilon_i)=\sigma^2$ which
are independent of
$W_{ij}$ and $Z_{ij}$.

By definition, $W_{ij}$ and $Z_{ij}$ possess substantially different
interpretations. $Z_{ij}$ describes the part of $X_{ij}$ which is
\textit{uncorrelated with all other variables}. A nonzero coefficient
$\beta
_j\neq0$
then means that the variation of $X_{ij}$ has a \textit{specific}
effect on
$Y_i$. We will of course assume that such nonzero coefficients are
\textit{sparse},
$\sharp\{j|\beta_j\neq0\}\leq S$ for some $S\ll p$. The
true variables $Z_{ij}$ are unknown, but with $\beta^{**}_j=\beta
_j^*-\beta_j$
model (\ref{eq-30}) can obviously be rewritten in the form
%
%
\begin{equation}\label{eq-31}
Y_i = \sum_{j=1}^p\beta^{**}_jW_{ij}+\sum_{j=1}^p\beta
_jX_{ij}+\varepsilon_i,\qquad
i=1,\ldots,n.
\end{equation}

The variables $W_{ij}$ are heavily correlated. It therefore does not
make any sense to assume that for some $j\in\{1,\ldots,p\}$ any
particular variable $W_{ij}$ possesses a specific influence on the
predictor variable. However,\vspace*{1pt} the term
$\sum_{j=1}^p\beta^{**}_jW_{ij}$ may represent an important,
\textit{common} effect of all predictor variables. The vectors $\mathbf{W}_i$
can obviously
be rewritten in terms of principal components.
Let us recall that $\lambda_1\geq\lambda_2\geq\cdots$ denote the
eigenvalues of the standardized covariance matrix of $\mathbf{W}_i$,
$\frac{1}{p}\bolds{\Gamma}=\frac{1}{p}\mathbb{E}(\mathbf
{W}_i\mathbf{W}_i^T)$ and $\bolds{\psi}_1,\bolds{\psi}_2,\ldots
$ corresponding orthonormal eigenvectors. We have
\[
\mathbf{W}_i = \sum_{r=1}^p(\bolds{\psi}^T_r\mathbf
{W}_i)\bolds{\psi}_r\quad
\mbox{and}\quad \sum_{j=1}^p\beta^{**}_jW_{ij}=\sum_{r=1}^p\alpha_{r}^*
(\bolds{\psi}^T_r\mathbf{W}_i),
\]
where $\alpha_r^*=\sum_{j=1}^p \beta_j^{**} \psi_{rj}$. As outlined in
the previous
sections we now
assume that the use of principal components allows for a considerable
reduction of
dimensionality, and that a small number of leading principal components
will suffice
to describe the effects of the variable $\mathbf{W}_i$. This may be
seen as an analogue of the sparseness assumption made for the $Z_{ij}$.
More precisely, subsequent analysis will be based on the assumption
that the
following \textit{augmented model} holds for some suitable $k\geq1$:
%
%
\begin{equation}\label{augmented-model}
Y_i = \sum_{r=1}^k\alpha_r\xi_{ir} + \sum_{j=1}^p\beta
_jX_{ij}+\varepsilon_i,
\end{equation}
where $\xi_{ir}=\bolds{\psi}^T_r\mathbf{W}_i/\sqrt{p\lambda_r}$ and
$\alpha_r=\sqrt{p\lambda_r}\alpha_r^*$. The use of $\xi_{ir}$ instead
of $\bolds{\psi}^T_r\mathbf{W}_i$ is motivated by the fact that
$\operatorname{Var}(\bolds{\psi}^T_r\mathbf{W}_i)=p\lambda_r$,
$r=1,\ldots,k$. Therefore the $\xi_{ir}$ are standardized variables
with $\operatorname{Var}(\xi_{i1})=\cdots
=\operatorname{Var}(\xi_{i1})=1$. Fitting an augmented model requires
us to select an appropriate $k$ as well as to determine sensible
estimates of $\xi_{i1},\ldots,\xi_{ik}$. Furthermore, model selection
procedures like Lasso or Dantzig have to be applied in order to
retrieve the nonzero coefficients $\alpha_r$, $r=1,\ldots,k$, and $\beta
_j$, $j=1,\ldots,p$. These issues will be addressed in subsequent
sections.

Obviously, the augmented model may be considered as a synthesis of the
standard type
of models proposed in the literature on functional regression and model
selection. It generalizes the classical multivariate linear
regression model (\ref{linear-model}). If a $k$-factor model holds
exactly, that is,
$\operatorname{rank}(\bolds{\Gamma})=k$, then the only substantial restriction
of (\ref{eq-30})--(\ref{augmented-model}) consists in the assumption
that $Y_i$ depends \textit{linearly} on
$W_i$ and $Z_i$.

We want to emphasize, however, that our analysis does not require the
validity of
a $k$-factor model. It is only assumed that there exists ``some'' $Z_i$
and $W_i$
satisfying our assumptions which lead
to (\ref{augmented-model}) for a sparse set of coefficients $\beta_j$.

\subsection{Identifiability}\label{sec31}

Let $\bolds{\beta}=(\beta_1,\ldots,\beta_p)^T$ and $\bolds
{\alpha}=(\alpha_1,\ldots,\alpha_k)^T$.
Since $\bolds{\psi}_r$, $r=1,\ldots,k$, are eigenvectors of
$\bolds{\Gamma}$ we have $\mathbb{E}(\bolds{\psi}^T_r\mathbf
{W}_i\bolds{\psi}_s^T\mathbf{W}_i)=0$ for
all $r,s=1,\ldots,p$, $r\ne s$. By assumption the random vectors
$\mathbf
{W}_i$ and $\mathbf{Z}_i$ are uncorrelated, and hence
$\mathbb{E}(\bolds{\psi}^T_r\mathbf{W}_iZ_{ij})=0$ for all
$r,j=1,\ldots,p$. Furthermore, $\mathbb{E}(Z_{il}Z_{ij})=0$ for all
$l\ne j$. If the augmented model (\ref{augmented-model}) holds, some
straightforward computations then show that under (A.1) for any alternative
set of coefficients $\bolds{\beta}^*=(\beta_1^*,\ldots,\beta
_p^*)^T$, $\bolds{\alpha}^*=(\alpha_1^*,\ldots,\alpha_k^*)^T$,
%
%
\begin{eqnarray}\label{abident}
&&\mathbb{E}\Biggl(\Biggl[\sum_{r=1}^k(\alpha_r-\alpha_r^*)\xi_{ir} +\sum
_{j=1}^p (\beta_j-\beta_j^*)X_{ij}
\Biggr]^2\Biggr)
\nonumber\\
&&\qquad\geq\sum_{r=1}^k\bigl(\alpha_r-\alpha_r^*+\sqrt{p\lambda_r}
\bolds{\psi}_r(\bolds{\beta}-\bolds{\beta}^*)\bigr)^2\\
&&\qquad\quad{}+D_1\| \bolds{\beta}-\bolds{\beta}^*\|^2_2.
\nonumber
\end{eqnarray}
We can conclude that the coefficients $\alpha_r$, $r=1,\ldots,k$ and
$\beta_j$, $j=1,\ldots,p$, in (\ref{augmented-model})
are uniquely determined.

Of course, an inherent difficulty of (\ref{augmented-model}) consists
of the fact that it contains the unobserved, ``latent'' variables
$\xi_{ir}=\bolds{\psi}^T_r\mathbf{W}_i/\sqrt{p\lambda_r}$.
To study this problem,
first recall that our setup imposes the decomposition (\ref
{var-decomp}) of the
covariance matrix $\bolds{\Sigma}$ of $\mathbf{X}_i$. If a factor
model with $k$ factors
holds exactly, then the $\bolds{\Gamma}$ possesses rank $k$. It
then follows from
well-established results in multivariate analysis that if $k<p/2$
the matrices $\bolds{\Gamma}$ and $\bolds{\Psi}$ are uniquely
identified.
If $\lambda_1>\lambda_2>\cdots>\lambda_k>0$, then also $\bolds{\psi
}_1,\ldots,\bolds{\psi}_k$ are
uniquely determined (up to sign) from the structure of $\bolds
{\Gamma}$.

However, for large $p$, identification is possible under even more
general conditions.
It is not necessary that a $k$-factor model holds exactly. We only
need an additional assumption on the magnitude
of the eigenvalues of $\frac{1}{p}\bolds{\Gamma}$ defining the
$k$ principal components of $\mathbf{W}_i$ to be considered.

\begin{longlist}[(A.3)]
\item[(A.3)] The eigenvalues of $\frac{1}{p}\bolds{\Gamma}$ are
such that
\[
\min_{j,l\leq k, j\neq l} |\lambda_j-\lambda_l|\ge v(k),\qquad
\min_{j\leq k} \lambda_j\ge v(k)
\]
for some $1\geq v(k)> 0$ with $pv(k)>6D_2$.
\end{longlist}

In the following we will qualitatively assume that $k\ll p$ as well as
$v(k)\gg1/p$. More specific assumptions will be made in the sequel.
Note that
eigenvectors are only unique up to sign
changes. In the following we will always assume that the right
``versions'' are used. This will go without saying.
\begin{theorem}\label{thmiden1}
Let $\xi_{ir}^*:=\frac{\bolds{\delta}^T_r\mathbf{X}_i}{\sqrt{p\mu _r}}$
and $\mathbf{P}_k = \mathbf{I}_p -
\sum_{j=1}^k\bolds{\psi}_j\bolds{\psi}_j^T$. Under assumptions
\textup{(A.1)} and \textup{(A.2)} we have\vadjust{\goodbreak} for all
$r=1,\ldots, k$, $j=1,\ldots,p$ and all $k,p$ satisfying
\textup{(A.3)}:
%
%
\begin{eqnarray}\label{thmiden-1}
|\mu_r-\lambda_r|&\leq&\frac{D_2}{p},
\\
\label{thmiden-2}
\|\bolds{\psi}_r-\bolds{\delta}_r\|_2&\leq&
\frac{2D_2}{pv(k)},\nonumber\\[-8pt]\\[-8pt]
\delta_{rj}^2&\leq&\frac{D_0}{pv(k)},\nonumber
\\
\label{thmiden-3}
\mathbb{E}([\xi_{ir}-
\xi_{ir}^*]^2)
&\leq&\frac{D_2}{p\mu_r}+\frac{(8\lambda_1+1)D_2^2}{p^2v(k)^2\mu_r},\nonumber\\[-8pt]\\[-8pt]
\mathbb{E}\biggl(\biggl[\xi_{ir}-
\frac{\bolds{\psi}^T_r\mathbf{X}_i}{\sqrt{p\lambda_r}}\biggr]^2\biggr)
&\leq&\frac{D_2}{p\lambda_r},\nonumber
\\
\label{thmiden-4}
\mathbb{E}\Biggl(\Biggl[\sum_{r=1}^k\alpha_r\xi_{ir} -
\sum_{r=1}^k\alpha_r\xi_{ir}^* \Biggr]^2\Biggr)
&\leq& k\sum_{r=1}^k\alpha_r^2\biggl(\frac{D_2}{p\mu_r}+
\frac{(8\lambda_1+1)D_2^2}{p^2v(k)^2\mu_r}\biggr).
\end{eqnarray}
\end{theorem}

For small $p$, standard factor analysis uses special algorithms in
order to identify~$\bolds{\psi}_r$. The theorem tells us that for large
$p$ this is unnecessary since then the eigenvectors $\bolds{\delta}_r$
of $\frac{1}{p}\bolds {\Sigma}$ provide a good approximation. The
predictor $\xi_{ir}^*$ of $\xi_{ir}$ possesses an error of order
$1/\sqrt{p}$. The error decreases as $p$ increases, and $\xi_{ir}^*$
thus yields a good approximation of $\xi_{ir}$ \textit{if $p$ is
large}. Indeed, if $p\rightarrow\infty$ [for fixed $\mu_r$, $v(k)$,
$D_1$ and $D_2$] then by (\ref{thmiden-3}) we have
$\mathbb{E}([\xi_{ir}-\xi _{ir}^*]^2)\rightarrow0$. Furthermore, by
(\ref{thmiden-4}) the error in predicting $\sum
_{r=1}^k\alpha_r\xi_{ir} + \sum_{j=1}^p\beta_jX_{ij}$ by $\sum
_{r=1}^k\alpha_r\xi_{ir}^* + \sum_{j=1}^p\beta_jX_{ij}$ converges to
zero as $p\rightarrow\infty$.

A crucial prerequisite for a reasonable analysis of the model is
sparseness of
the coefficients $\beta_j$. Note that if $p$ is large compared to $n$,
then by
(\ref{thmiden-4})
the error in replacing $\xi_{ir}$ by $\xi_{ir}^*$ is negligible
compared to the estimation error induced by the
existence of the error terms $\varepsilon_i$. If $k\ll p$ and
$\sharp\{j|\beta_j\neq0\}\ll p$, then the true coefficients $\alpha
_r$ and
$\beta_j$ provide a \textit{sparse} solution of the regression problem.

Established
theoretical results [see \citet{BiRiTss09}] show that under some
regularity conditions (validity of the ``restricted eigenvalue
conditions'') model selection procedures \textit{allow to identify such
sparse solutions} even if there are multiple vectors of coefficients
satisfying the normal equations. The latter is of course always the
case if $p>n$. Indeed, we will show in the following sections that
factors can be consistently estimated from the data, and that
a suitable application of Lasso or the Dantzig-selector leads to consistent
estimators $\widehat\alpha_r$, $\widehat\beta_j$ satisfying
$\sup_r |\alpha_r-\widehat\alpha_r|\rightarrow_P 0$, $\sup_j |\beta
_j-\widehat
\beta_j|\rightarrow_P 0$, as $n,p\rightarrow\infty$.

When replacing $\xi_{ir}$ by $\xi_{ir}^*$, there are alternative sets
of coefficients leading to the same prediction error as in
(\ref{thmiden-4}). This is due\vadjust{\goodbreak} to the fact that $\xi_{ir}^*=
\sum_{j=1}^p\frac{\delta_{rj}X_{ij}}{\sqrt{p\mu_r}}$. However, all
these alternative solutions are \textit{nonsparse} and cannot be
identified by Lasso or other procedures. In particular, it is easily
seen that
%
%
\begin{eqnarray}\label{stanmodident}
\sum_{r=1}^k\alpha_r\xi_{ir}^* + \sum_{j=1}^p\beta_jX_{ij}
=\sum_{j=1}^p\beta_j^{\mathit{LR}}X_{ij}\nonumber\\[-8pt]\\[-8pt]
&&\eqntext{\displaystyle \mbox{with }
\beta^{\mathit{LR}}_j:=\beta_j+ \sum_{r=1}^k\alpha_r\frac{\delta_{rj}}{\sqrt{p\mu_r}}.}
\end{eqnarray}

By (\ref{thmiden-2}) all values $\delta_{rj}^2$ are of order
$1/(pv(k))$. Since
$\sum_j \delta_{rj}^2=1$, this implies that many $\delta_{rj}^2$ are
nonzero. Therefore, if $\alpha_r\neq0$ for some $r\in\{1,\ldots,k\}$, then
$\{j|\beta_j^{\mathit{LR}}\neq0\}$ contains a large number of small,
nonzero coefficients
and is not at all sparse. If $p$ is large compared to $n$ no known
estimation procedure will be able to provide consistent estimates of
these coefficients.


Summarizing the above discussion we can conclude:

\begin{longlist}[(2)]
\item[(1)]
If the variables $X_{ij}$ are heavily correlated and follow an
approximate
factor model, then one may reasonably expect substantial effects of the
common, joint variation of all variables and, consequently, nonzero coefficients
$\beta_j^*$ and $\alpha_r$ in (\ref{eq-30}) and (\ref{augmented-model}).
But then a ``bet on sparsity'' is unjustifiable when dealing with the standard
regression model (\ref{linear-model}). It follows from
(\ref{stanmodident}) that for large $p$ model (\ref{linear-model})
holds approximately for a nonsparse set of coefficients $\beta
_j^{\mathit{LR}}$, since many
small, nonzero coefficients are necessary in order to capture the
effects of
the common joint variation.

\item[(2)]
The augmented model offers a remedy to this problem by
pooling possible
effects of the joint variation using a small number of additional variables.
Together with the familiar assumption of a small number of variables possessing
a specific influence, this leads to a sparse model with at most $k+S$ nonzero
coefficients which can be recovered from model selection procedures
like Lasso or
the Dantzig-selector.

\item[(3)]
In practice, even if
(\ref{augmented-model}) only holds approximately, since a too-small value
of $k$ has been selected, it may be able to quantify at least some
important part of the
effects discussed above. Compared to an analysis based on a standard
model~(\ref{linear-model}),
this may lead to a substantial improvement of model fit as well as to more
reliable interpretations of significant variables.
\end{longlist}

\subsection{Estimation}\label{sec32}
For a pre-specified $k\geq1$ we now define a procedure for estimating
the components of the corresponding augmented model (\ref{augmented-model})
from given data.
This obviously specifies suitable procedures
for approximating the unknown values $\xi_{ir}$ as well as to apply
subsequent model selection procedures in order to retrieve nonzero coefficients
$\alpha_r$ and $\beta_j$, $r=1,\ldots,k$, $j=1,\ldots,p$. A discussion
of the
choice of $k$ can be found in the next section.

Recall from Theorem \ref{thmiden1} that for large $p$ the eigenvectors
$\bolds{\psi}_1,\ldots,\bolds{\psi}_k$ of $\frac
{1}{p}\bolds{\Gamma}$
are well approximated by the
eigenvectors of the standardized covariance matrix $\frac
{1}{p}\bolds{\Sigma}$. This motivates us to use the
empirical principal components of $\mathbf{X}_1,\ldots,\mathbf{X}_n$ in order
to determine estimates of $\bolds{\psi}_r$ and $\xi_{ir}$.
Theoretical support
will be given in the next section.
Define\vspace*{1pt} $\widehat{\lambda}_1\geq\widehat{\lambda}_2\geq\cdots$ as the
eigenvalues of the standardized empirical covariance\vspace*{1pt} matrix $\frac
{1}{p}\widehat{\bolds{\Sigma}}=\frac{1}{np}\sum_{i=1}^n\mathbf
{X}_i^T\mathbf{X}_i$, while $\widehat{\bolds{\psi}}_1,\widehat
{\bolds{\psi}}_2,\ldots$ are associated orthonormal eigenvectors.
We then estimate $\xi_{ir}$ by
\[
\widehat{\xi}_{ir}=\widehat{\bolds{\psi}}_{r}^T
\mathbf{X}_i/\sqrt{p\widehat{\lambda}_r},\qquad r=1,\ldots,k, i=1,\ldots,n.
\]

When replacing $\xi_{ir}$ by $\widehat{\xi}_{ir}$ in
(\ref{augmented-model}), a direct application of model selection
procedures does not seem to be adequate, since $\widehat{\xi}_{ir}$ and
the predictor variables $X_{ij}$ are heavily correlated. We therefore
rely on a projected model. Consider the projection matrix on the
orthogonal space of the space spanned\vspace*{1pt} by the eigenvectors
corresponding to the $k$ largest eigenvalues of
$\frac{1}{p}\widehat{\bolds{\Sigma}}$
\[
\widehat{\mathbf{P}}_k = \mathbf{I}_p - \sum_{r=1}^k\widehat
{\bolds{\psi}}_r\widehat{\bolds{\psi}}{}^T_r.
\]
Then model (\ref{augmented-model}) can be rewritten for $i=1,\ldots,n$,
%
%
\begin{equation}\label{model-projected}
Y_i = \sum_{r=1}^k\widetilde{\alpha}_r\widehat{\xi}_{ir} + \sum
_{j=1}^p\widetilde{\beta}_j\frac{(\widehat{\mathbf{P}}_k\mathbf
{X}_i)_j}{(({1/n})\sum_{i=1}^n(\widehat{\mathbf{P}}_k\mathbf
{X}_i)_j^2)^{1/2}} + \widetilde{\varepsilon}_i+\varepsilon_i,
\end{equation}
where $\widetilde{\alpha}_r=\alpha_r+\sqrt{p\widehat{\lambda}_r}
\sum_{j=1}^p\widehat{\psi}_{rj}\beta_j$, $\widetilde{\beta}_j=$
$\beta_j(\frac{1}{n}\sum_{i=1}^n(\widehat{\mathbf{P}}_k\mathbf
{X}_i)_j^2)^{1/2}$ and $\widetilde{\varepsilon}_i=\sum_{r=1}^k\alpha
_r(\xi_{ir}-\widehat{\xi}_{ir})$.
It will be shown in the next section that for large $n$ and $p$ the additional
error term $\widetilde\varepsilon$ can be assumed to be reasonably small.

In the following we will use $\widetilde{\mathbf{X}}_i$ to denote the
vectors with entries $\widetilde{X}_{ij}:=
\frac{(\widehat{\mathbf{P}}_k\mathbf{X}_i)_j}{( ({1/n})\sum
_{l=1}^n(\widehat{\mathbf{P}}_k\mathbf{X}_l)^2_j)^{1/2}}$.
Furthermore,\vspace*{1pt} consider the $(k+p)$-dimensional vector of
predictors
$\bolds{\Phi}_i:=(\widehat{\xi}_{i1},\ldots,\widehat{\xi}_{ik}$,
$\widetilde{X}_{i1},\ldots,\widetilde{X}_{ip})^T$. The Gram matrix in
model (\ref{model-projected}) is a block matrix defined as
\[
\frac{1}{n}\sum_{i=1}^n\bolds{\Phi}_i\bolds{\Phi}_i^T = \pmatrix{
\mathbf{I}_k&\bolds{0}\vspace*{2pt}\cr \bolds{0}&\displaystyle
\frac{1}{n}\sum_{i=1}^n\widetilde{\mathbf{X}}_i\widetilde
{\mathbf{X}}_i^T},
\]
where $\mathbf{I}_k$ is the identity matrix of size $k$. Note that the
normalization of the predictors in (\ref{model-projected}) implies that
the diagonal elements of the Gram matrix above are equal to 1.

Arguing now that the vector of parameters $\bolds{\theta
}:=(\widetilde{\alpha}_1,\ldots,\widetilde{\alpha}_k,\widetilde{\beta
}_1,\ldots,\widetilde{\beta}_p)^T$ in model (\ref{model-projected}) is
$(k+S)$-sparse, we may use a selection procedure to recover/estimate
the nonnull parameters. In the following we will concentrate on the
Lasso estimator introduced in \citet{Ti96}. For a pre-specified
parameter $\rho>0$, an estimator $\bolds{\theta}$ is then obtained as
%
%
\begin{equation}\label{dantzig2}
\widehat{\bolds{\theta}} = \mathop{\arg\min}_{\widetilde{\bolds
{\theta}}\in\mathbb{R}^{k+p}}
\frac{1}{n}\|\mathbf{Y}-\bolds{\Phi}\widetilde{\bolds{\theta
}}\|_2+2\rho\|\widetilde{\bolds{\theta}}\|_1,
\end{equation}
$\bolds{\Phi}$ being the $n\times(k+p)$-dimensional matrix with
rows $\bolds{\Phi}_i$. We can alternatively use the Dantzig
selector introduced in \citet{CaTa07}.

Finally, from $\widehat{\bolds{\theta}}$, we define corresponding
estimators for $\alpha_r$, $r=1,\ldots,k$, and $\beta_j$,
$j=1,\ldots,p$,
in the unprojected model (\ref{augmented-model}).
\[
\widehat{\beta}_j =
\frac{\hspace*{1.5pt}\widehat{\hspace*{-1.5pt}\widetilde{\beta}}_j}{(
({1/n})\sum_{i=1}^n(\widehat{\mathbf{P}}_k\mathbf{X}_i)_j^2
)^{1/2}},\qquad
j=1,\ldots,p,
\]
and
\[
\widehat{\alpha}_r = \widehat{\widetilde{\alpha}}_r-\sqrt{p\widehat
{\lambda}_r}\sum_{j=1}^p\widehat{\psi}_{rj}\widehat{\beta}_j,\qquad
r=1,\ldots,k.
\]

\section{High-dimensional factor analysis: Theoretical results}\label{sec4}

The following theorem shows that principal components which are able
to explain a considerable proportion of total variance can be estimated
consistently.

For simplicity, we will concentrate on the case
that $n$ as well as $p>\sqrt{n}$ are large enough such that
\begin{longlist}[(A.4)]
\item[(A.4)] $C_0 (\log p /n)^{1/2}\ge\frac{D_0}{p}$ and
$v(k)\geq6 (D_2/p+C_0(\log p/n)^{1/2})$.
\end{longlist}
\begin{theorem}\label{thmapp1}
Under assumptions \textup{(A.1)--(A.4)} and under events
(\ref{B21})--(\ref {B24}) we have for all $r=1,\ldots, k$ and all
$j=1,\ldots,p$,
%
%
\begin{eqnarray}
\label{thmapp1-1}
|\lambda_r-\widehat{\lambda}_r|&\leq&\frac{D_2}{p}+C_0(\log p/n)^{1/2}
\\
\label{thmapp1-2}
\|\bolds{\psi}_r-\widehat{\bolds{\psi}}_r\|_2&\leq&
2\frac{{D_2}/{p}+C_0(\log p/n)^{1/2}}{v(k)},
\\
\label{thmapp1-3}
\psi_{rj}^2&\leq&\frac{D_0-D_1}{p\lambda_r}\leq\frac
{D_0-D_1}{pv(k)},\\
\label{thmapp1-4}
\widehat{\psi}_{rj}^2&\leq&\frac{D_0+C_0(\log p/n)^{1/2}}{p\widehat
{\lambda}_r}\nonumber\\[-8pt]\\[-8pt]
&\leq&\frac{6}{5} \frac{D_0+C_0(\log p/n)^{1/2}}{pv(k)}.\nonumber
\end{eqnarray}
\end{theorem}

Theorem \ref{thmapp1} shows that for sufficiently large $p$ ($p>\sqrt
{n}$) the eigenvalues and eigenvectors of $\frac{1}{p}\widehat
{\bolds{\Sigma}}$
provide reasonable estimates of $\lambda_r$ and $\bolds{\psi}_r$
for $r=1,\ldots,k$.
Quite obviously it is not possible to determine sensible estimates of
\textit{all}
$p$ principal components of $\frac{1}{p}\bolds{\Sigma}$. Following
the proposition it is required that $\lambda_r$ as well as $\mu_r$ be
of order at least
$\sqrt{\frac{\log p}{n}}$. Any smaller component cannot be
distinguished from
pure ``noise'' components. Up to the $\log p$-term this corresponds to
the results
of \citet{HaHos07} who study the problem of the
number of principal components that can be consistently estimated in a
functional
principal component analysis.

The above insights are helpful for selecting an appropriate $k$ in a
real data
application. In tendency, a suitable factor model will incorporate $k$
components which explain a large percentage of the total variance of
$\mathbf{X}_i$, while $\lambda_{k+1}$ is very small.
If for a sample of high-dimensional vectors $\mathbf{X}_i$ a principal
component analysis leads to the conclusion that the first (or second,
third$,\ldots$) principal components explains a \textit{large}
percentage of
the total (empirical) variance of the observations, then such a
component cannot be generated by ``noise'' but reflects an underlying
structure. In particular, such a component
may play a crucial role in modeling a response variable $Y_i$ according to
an augmented regression model of the form
(\ref{augmented-model}).

\citet{BaiNg02} develop criteria of selecting the dimension $k$ in a
high-dimensional factor model. They rely on an adaptation of the well-known
AIC and BIC procedures in model selection. One possible approach is as
follows: Select a maximal possible dimension $k_{\max}$ and estimate
$\bar\sigma^2=
\frac{1}{p}\sum_{j=1}^p \sigma_j^2$ by $\widehat\sigma^2=
\frac{1}{np}\sum_{i=1}^n\sum_{j=1}^p (X_{ij}-\sum_{r=1}^{k_{\max}}
(\widehat{\bolds{\psi}}{}^{T}_{r}\mathbf{X}_i)\widehat\psi_{rj})^2$.
Then\vspace*{1pt} determine an estimate $\widehat k$ by minimizing
%
%
\begin{equation}\label{baicrit}
\frac{1}{np}\sum_{i=1}^n\sum_{j=1}^p \Biggl(X_{ij}-\sum_{r=1}^{\kappa}
(\widehat{\bolds{\psi}}_{r}\mathbf{X}_i)\widehat\psi_{rj}\Biggr)^2
+\kappa\widehat\sigma^2 \biggl( \frac{n+p}{np}\biggr)
\log\min\{n,p\}
\end{equation}
over $\kappa=1,\ldots,k_{\max}$. \citet{BaiNg02} show that
under some regularity conditions this criterium (as well as a number of
alternative versions) provides asymptotically consistent estimates
of the true factor dimension $k$ as $n,p\rightarrow\infty$. In our
context these
regularity conditions are satisfied if (A.1)--(A.4) hold for all $n$
and $p$,
$\sup_{j,p} \mathbb{E}(Z_{ij}^8)<\infty$ and if there exists some
$B_0>0$ such
that $\lambda_k\geq B_0>0$,
for all $n,p$.

Now recall the modified version (\ref{model-projected}) of the
augmented model used
in our estimation procedure. The following theorem establishes bounds for
the projections $(\widehat{\mathbf{P}}_k\mathbf{X}_i)_j$ as well as
for the
additional error terms $\widetilde{\varepsilon}_i$.
Let $\mathbf{P}_k = \mathbf{I}_p -
\sum_{j=1}^k\bolds{\psi}_j\bolds{\psi}_j^T$ denote the
population version of $\widehat{\mathbf{P}}_k$.
\begin{theorem}\label{thmapp2}
Assume \textup{(A.1)} and \textup{(A.2)}. There then exist constants
$M_1$, $ M_2$, $M_3<\infty$, such that for all $n,p,k$ satisfying
\textup{(A.3)} and \textup{(A.4)}, all $j,l\in\{1,\ldots,p\}$, $j\neq
l$,
%
%
\begin{eqnarray}
\label{pkapp0}
\frac{1}{n}\sum_{i=1}^n (\widehat{\mathbf{P}}_k\mathbf{X}_i)_j^2
&\geq&
\sigma_j^2-
M_1\frac{kn^{-1/2}\sqrt{\log p}}{v(k)^{1/2}} ,
\\
\label{pkapp1}
\Biggl|\frac{1}{n}\sum_{i=1}^n (\widehat{\mathbf{P}}_k\mathbf
{X}_i)_j^2-\sigma_j^2\Biggr|&\leq&\mathbb{E}((\mathbf{P}_k\mathbf
{W}_i)_j^2 )+
M_2\frac{kn^{-1/2}\sqrt{\log p}}{v(k)^{3/2}}
\end{eqnarray}
hold with probability $A(n,p)$, while 
%
%
\begin{eqnarray} \label{pkapp2}
\frac{1}{n}\sum_{i=1}^n \widetilde{\varepsilon}_i^2&=&
\frac{1}{n}\sum_{i=1}^n\Biggl(\sum_{r=1}^k(\widehat{\xi}_{ir}-\xi
_{ir})\alpha_r\Biggr)^2\nonumber\\[-8pt]\\[-8pt]
&\leq&\frac{k\alpha^2_{\mathrm{sum}}M_3}{v(k)^3}\biggl(\frac{\log p}{n}+
\frac{v(k)^2}{p}\biggr)\nonumber
\end{eqnarray}
holds with probability at least $A(n,p)-\frac{k}{n}$. Here,
$\alpha^2_{\mathrm{sum}}=\sum_{r=1}^k \alpha_r^2$.
\end{theorem}

Note that if $\mathbf{X}_i$ satisfies a $k$-dimensional factor model,
that is, if the rank of
$\frac{1}{p}\bolds{\Gamma}$ is equal to $k$, then $\mathbf
{P}_k\mathbf{W}_i=0$. The theorem then states that for large $n$ and
$p$ the projected variables $(\widehat{\mathbf{P}}_k\mathbf{X}_i)_j$,
$j=1,\ldots,k$,
``in average'' behave similarly to the specific variables $Z_{ij}$.
Variances will be close
to $\sigma_j^2=\operatorname{Var}(Z_{ij})$.

\section{Theoretical properties of the augmented model}\label{sec5}

We come back to model~(\ref{augmented-model}). As shown in Section \ref{sec32},
the Lasso or the Dantzig selector may be used to determine estimators
of the parameters of the model. Identification of sparse solutions as
well as
consistency of estimators require structural assumptions on the
explanatory variables. The weakest assumption on the correlations
between different variables seems to be the so-called \textit
{restricted eigenvalue} condition introduced by \citet{BiRiTss09};
see also \citet{ZhVaBu09}.

We first provide
a theoretical result which shows that for large $n,p$ the design
matrix of the
projected model
(\ref{model-projected}) satisfies the restricted eigenvalue conditions
given in \citet{BiRiTss09} with high probability. We will
additionally assume that $n,p$ are large enough such that
\begin{longlist}[(A.5)]
\item[(A.5)] $D_1/2>M_1\frac{kn^{-1/2}\sqrt{\log p}}{v(k)^{1/2}}$,
\end{longlist}
where $M_1$ is defined as in Theorem \ref{thmapp2}.

Let $J_0$ denote an arbitrary subset of indices, $J_0\subset\{1,\ldots
,p\}$ with $|J_0|\le k+S$. For a vector $\mathbf{a}\in\mathbb
{R}^{k+p}$, let $\mathbf{a}_{J_0}$ be the vector in $\mathbb{R}^{k+p}$
which has the same coordinates as $\mathbf{a}$ on $J_0$ and zero
coordinates on the complement $J_0^c$ of $J_0$. We define in the same
way $\mathbf{a}_{J_0^c}$. Now for $k+S\le(k+p)/2$ and for an integer
$m\ge k+S$, $S+m\le p$, denote by $J_m$ the subset of $\{1,\ldots,k+p\}
$ corresponding to $m$ largest in absolute value coordinates of $\mathbf
{a}$ outside of $J_0$, and define $J_{0,m}:=J_0\cup J_m$. Furthermore,
let $(x)_+=\max\{x,0\}$.
\begin{prop} \label{propsparse2}
Assume \textup{(A.1)} and \textup{(A.2)}. There then exists
a constant $M_4<\infty$, such that for all
$n,p,k,S$, $k+S\le(k+p)/2$, satisfying \textup{(A.3)--(A.5)},
and $c_0=1,3$
%
%
\begin{eqnarray}\label{sparse02}
&&\kappa(k+S,k+S,c_0)\nonumber\\[-2pt]
&&\quad:=\min_{J_0\subset\{1,\ldots,k+p\}\dvtx|J_0|\leq k+S}
\min_{\bolds{\Delta}\neq0\dvtx\|\bolds{\Delta}_{J_0^c}\|_1\leq c_0\|
\bolds{\Delta}_{J_0}\|_1}
\frac{[\bolds{\Delta}^T({1}/{n})\sum_{i=1}^n \bolds{\Phi
}_i\bolds{\Phi}_i^{T} \bolds{\Delta}]^{1/2}}{\|\bolds{\Delta
}_{J_{0,k+S}}\|_2}\hspace*{-12pt}
\nonumber\\[-8pt]\\[-8pt]
&&\quad\geq\biggl(\frac{D_1}{D_0+C_0n^{-1/2}\sqrt{\log p}}
-\frac{8(k+S)c_0M_4k^2n^{-1/2}\sqrt{\log
p}}{v(k)(D_1-M_1kv(k)^{1/2}n^{-1/2}\sqrt{\log p})}\biggr)_+^{{1}/{2}}
\nonumber\\[-2pt]
&&\quad
=:\mathbf{K}_{n,p}(k,S,c_0)\nonumber
\end{eqnarray}
holds with probability $A(n,p)$.
\end{prop}

Asymptotically, if $n$ and $p$ are large, then $\mathbf
{K}_{n,p}(k,S,c_0)>0$, $c_0=1,3$, provided that
$k$, $S$ and $1/v(k)$ are sufficiently small compared to $n,p$. In this
case the proposition implies
that with high probability the restricted eigenvalue condition
$\operatorname{RE}(k+S,k+S,c_0)$ of \citet{BiRiTss09} [i.e., $\kappa
(k+S,k+S,c_0)>0$] is satisfied. The same holds for the conditions
$\operatorname{RE}(k+S,c_0)$ which
require $\kappa(k+S,c_0)>0$, where\vadjust{\eject}
\begin{eqnarray*}
&&
\kappa(k+S,c_0)\\[-3pt]
&&\qquad:=\min_{J_0\subset\{1,\ldots,k+p\}\dvtx|J_0|\leq k+S}
\min_{\bolds{\Delta}\neq0\dvtx\|\bolds{\Delta}_{J_0^c}\|_1\leq c_0\|
\bolds{\Delta}_{J_0}\|_1}
\frac{[\bolds{\Delta}^T({1/n})\sum_{i=1}^n \bolds{\Phi
}_i\bolds{\Phi}_i^{T} \bolds{\Delta}]^{1/2}}{\|\bolds{\Delta
}_{J_{0}}\|_2}\\[-3pt]
&&\qquad\geq\kappa(k+S,k+S,c_0).
\end{eqnarray*}

The following theorem now provides bounds for the $L^1$ estimation
error and the $L^2$ prediction loss for the Lasso estimator of the
coefficients of the augmented model. It generalizes the results of
Theorem 7.2 of \citet{BiRiTss09} obtained under the standard
linear regression model. In our analysis merely the values of
$\kappa(k+S,c_0)$ for $c_0=3$ are of interest. However, only slight
adaptations of the proofs are necessary in order to derive
generalizations of the bounds provided by \citet{BiRiTss09} for
the Dantzig selector ($c_0=1$) and for the $L^q$ loss, $1<q\le2$. In
the latter case, $\kappa(k+S,c_0)$ has to be replaced by
$\kappa(k+S,k+S,c_0)$. In the following, let $M_1$ and $M_3$ be defined
as in Theorem~\ref{thmapp2}.
\begin{theorem}\label{thmaug} Assume \textup{(A.1)}, \textup{(A.2)} and
suppose that the error terms
$\varepsilon_i$ in model (\ref{augmented-model})
are independent $\mathcal{N}(0,\sigma^2)$ random variables with $\sigma
^2>0$. Now consider the Lasso estimator $\widehat{\bolds{\theta}}$
defined by (\ref{dantzig2}) with
\[
\rho= A\sigma\sqrt{\frac{\log(k+p)}{n}}+\frac{M_5\alpha
_{\mathrm{sum}}}{v(k)^{3/2} }\sqrt{\frac{\log p}{n}},
\]
where
$A>2\sqrt{2}$, $M_5$ is a positive constant and $\alpha_{\mathrm{sum}}=\sum
_{r=1}^k|\alpha_r|$.

\enlargethispage{16pt}
If $M_5<\infty$ is sufficiently large, then for all
$n,p,k$, $k+S\le(k+p)/2$, satisfying \textup{(A.3)--(A.5)} as well
as $\mathbf{K}_{n,p}(k,S,3)>0$, the following inequalities hold with
probability at least $A(n,p)-(p+k)^{-A^2/2}$:
%
%
\begin{eqnarray}
\label{thmaug-1}\quad
\sum_{r=1}^k|\widehat{\alpha}_r-\alpha_r| &\le& \frac{16(k+S)}{\kappa
^2}\nonumber\\[-8pt]\\[-8pt]
&&{}\times\rho\biggl(1+\frac{k(D_0+C_0n^{-1/2}\sqrt{\log
p})^{1/2}}{(D_1-M_1({kn^{-1/2}\sqrt{\log
p}}/{v(k)^{1/2}}))^{1/2}}\biggr),\nonumber\\
\label{thmaug-2}
\sum_{j=1}^p|\widehat{\beta}_j-\beta_j| &\le& \frac{16(k+S)}{\kappa
^2(D_1-M_1({kn^{-1/2}\sqrt{\log p}}/{v(k)^{1/2}}))^{1/2}}\rho,
\end{eqnarray}
where $\kappa=\kappa(k+S,3)$. Moreover,
%
%
\begin{eqnarray}\label{thmaug2-1}
&&\frac{1}{n}\sum_{i=1}^n\Biggl(\sum_{r=1}^k\widehat{\xi}_{ir}\widehat
{\alpha}_r
+\sum_{j=1}^pX_{ij}\widehat{\beta}_j- \Biggl(\sum_{r=1}^k\xi_{ir}\alpha_r+\sum
_{j=1}^pX_{ij}\beta_j\Biggr)\Biggr)^2\nonumber\\[-9pt]\\[-9pt]
&&\qquad\le\frac{32(k+S)}
{\kappa^2}\rho^2+\frac{2k\alpha^2_{\mathrm{sum}}M_3}{v(k)^3}\biggl(\frac{\log p}{n}+
\frac{v(k)^2}{p}\biggr)\nonumber
\end{eqnarray}
holds with probability at least $A(n,p)-(p+k)^{-A^2/2}-\frac{k}{n}$.
\end{theorem}
\eject

Of course, the main message of the theorem is asymptotic in nature. If
$n,p$ tend to infinity
for fixed values of $k$ and $S$, then the $L_1$ estimation error and
the $L^2$ prediction error converge at rates
$\sqrt{\log p/n}$
and $\log p/n+1/p$, respectively.
For values of $k$ and $S$ tending to infinity as the sample size tends
to infinity, the rates are more complicated. In particular, they depend
on how fast $v(k)$ converges to
zero as $k\rightarrow\infty$.
Similar results hold for the estimators based on the Dantzig selector.
\begin{rem}
Note that Proposition \ref{propsparse2} as well as the results of
Theorem \ref{thmaug} heavily depend on the validity of assumption (A.1)
and the corresponding value $0<D_1\leq\inf_j \sigma_j^2$, where
$\sigma_j^2=\operatorname{var}(Z_{ij})$. It is\vspace*{1pt} immediately
seen that the smaller the $D_1$ the smaller the value of
$\kappa(k+S,k+S,c_0)$ in (\ref{sparse02}). This means that \textit{all}
variables $X_{ij}=W_{ij}+Z_{ij}$, $j=1,\ldots,p$ have to possess a
sufficiently large specific variation which is not shared by other
variables. For large $p$ this may be seen as a restrictive assumption.
In such a situation one may consider a restricted version of model
(\ref {augmented-model}), where variables with extremely small values
of $\sigma_j^2$ are eliminated. But for large $n,p$ we can infer from
Theorem \ref{thmapp2} that a small value of
$\frac{1}{n}\sum_{i=1}^n(\widehat{\mathbf{P}}_k\mathbf {X}_i)^2_j$
indicates that also $\sigma_j^2$ is small. Hence, an extension of our
method consists of introducing some threshold $D_{\mathrm{thresh}}>0$
and discarding\vadjust{\goodbreak} all those variables $X_{ij}$,
$j\in\{1,\ldots,p\}$, with $\frac{1}{n}\sum_{i=1}^n(\widehat{\mathbf
{P}}_k\mathbf{X}_i)^2_j<D_{\mathrm{thresh}}$. A precise analysis is not
in the scope of the present paper.
\end{rem}
\begin{rem}
If $\alpha_1=\cdots=\alpha_k=0$
the augmented model reduces to the
standard linear regression model (\ref{linear-model}) with a
sparse set
of coefficients, $\sharp\{j|\beta_j\neq0\}\le S$ for some $S\le p$.
An application of our estimation procedure is then unnecessary, and coefficients
may be estimated by traditional model selection procedures. Bounds on
estimation errors can therefore be directly obtained from the results
of \citet{BiRiTss09}, provided that the restricted eigenvalue
conditions are satisfied.
But
in this situation
a slight adaptation of the proof of Proposition \ref{propsparse2}
allows us to establish
a result similar to (\ref{sparse02}) for the standardized variables
$X^*_{ij}:=X_{ij}/(\frac{1}{n}\sum_{i=1}^nX_{ij}^2)^{1/2}$.
Define $\mathbf{X}^*$ as the $n\times p$-matrix with generic elements
$X_{ij}^*$. When\vspace*{1pt} assuming (A.1), (A.2) as well as $D_1-3C_0n^{-1/2}\sqrt
{\log p}>0$, then for $S\le p/2$ the following inequality holds with
probability $A(n,p)$:
%
%
\begin{eqnarray}\label{sparse01}\quad
&&\kappa(S,S,c_0)\nonumber\\
&&\qquad:=\min_{J_0\subset\{1,\ldots,p\}\dvtx|J_0|\leq S}
\min_{\bolds{\Delta}\neq0\dvtx\|\bolds{\Delta}_{J_0^c}\|_1\leq c_0\|
\bolds{\Delta}_{J_0}\|_1}
\frac{[\bolds{\Delta}^T({1/n})\sum_{i=1}^n \mathbf{X}^*_i\mathbf
{X}_i^{*T} \bolds{\Delta}]^{1/2}}{\|\bolds{\Delta}_{J_{0,S}}\|_2}
\\
&&\qquad\geq\biggl(\frac{D_1}{D_0+C_0n^{-1/2}\sqrt{\log p}}-\frac
{8Sc_0C_0 n^{-1/2}\sqrt{\log p}}{D_1-3C_0n^{-1/2}\sqrt{\log p}}\biggr)
_+^{1/2},
\nonumber
\end{eqnarray}
where $c_0=1,3$. Recall, however, from the discussion in Section \ref{sec31}
that $\alpha_1=\cdots=\alpha_k=0$ is a restrictive condition in the
context of highly correlated regressors.
\end{rem}

\section{Simulation study}\label{sec6}

In this section we study the finite sample performance of the
estimators discussed in the proceeding sections. We consider a factor
model with $k=2$ factors. The first
factor is
$\psi_{1j}=1/\sqrt{p}$, $j=1,\ldots,p$, while the second factor is
given by
$\psi_{2j}=1/\sqrt{p}$, $j=1,\ldots,p/2$, and $\psi_{2j}=-1/\sqrt{p}$,
$j=p/2+1,\ldots,p$. For different values of $n,p,\alpha_1,\alpha_2$ and
$0<\lambda_1<1, 0<\lambda_2<1$ observations $(\mathbf{X}_i,Y_i)$ with
$\operatorname{var}(X_{ij})=1$, $j=1,\ldots,p$,
are generated according to the
model
%
%
\begin{eqnarray}
\label{simul-1}
X_{ij}&=&\sqrt{p\lambda_1}\xi_{i1}\psi_{1j}+\sqrt{p\lambda_2}\xi
_{i2}\psi
_{2j}+Z_{ij},
\\
\label{simul-2}
Y_i&=&\alpha_1 \xi_{i1}+\alpha_2 \xi_{i1}+\sum_{j=1}^p \beta_j
X_{ij}+\varepsilon_i,
\end{eqnarray}
where $\xi_{ir}\sim N(0,1)$, $r=1,2$, $Z_{ij}\sim N(0,1-\lambda
_1-\lambda_2)$, and
$\varepsilon_{i}\sim N(0,\sigma^2)$ are independent variables. Our study
is based
on $S=\sharp\{j|\beta_j\neq0\}=4$ nonzero $\beta$-coefficients whose values
are $\beta_{10}=1$, $\beta_{20}=0.3$, $\beta_{21}=-0.3$ and $\beta
_{40}=-1$, while
the error variance is set to $\sigma^2=0.1$.

The parameters of the augmented model with $k=2$ are estimated by
using the Lasso-based estimation procedure described in Section \ref{sec32}.
The behavior of the estimates is compared to the Lasso estimates of the
coefficients of a standard regression model (\ref{linear-model}).
All results reported in this section are obtained by applying the
LARS-package by Hastie and Efron implemented in R. All tables
are based on 1,000 repetitions of the simulation experiments.
The corresponding R-code can be obtained from the authors upon request.

\begin{figure}

\includegraphics{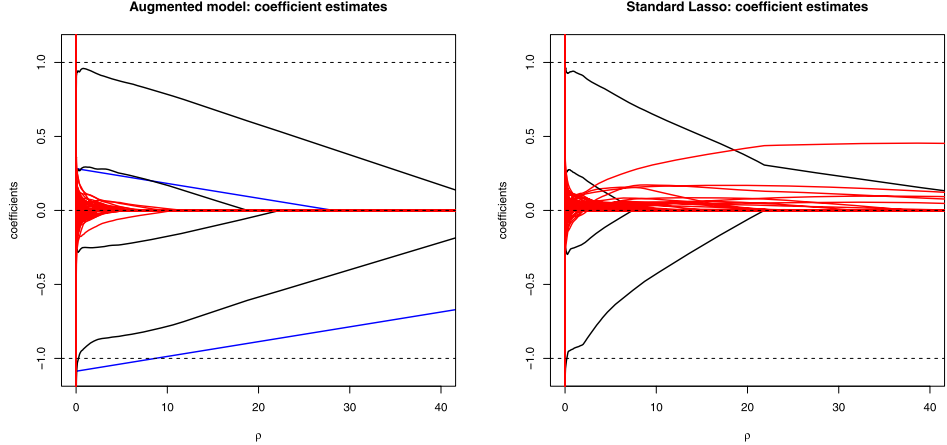}

\caption{Paths of Lasso estimates for the augmented model (left panel) and
the standard linear model (right panel) in dependence of $\rho$;
black---estimates
of nonzero $\beta_j$; red---estimates of coefficients with $\beta_j=0$;
blue---$\widehat{\widetilde{\alpha}}_r$.} \label{figur1}
\end{figure}

Figure \ref{figur1} and Table \ref{tabLa1} refer to the situation with $\lambda_1=0.4$,
$\lambda_2=0.2$, $\alpha_1=1$ and $\alpha_2=-0.5$. We then
have $\operatorname{var}(X_{ij})=1$, $j=1,\ldots,p$, while the first and second factor
explain 40\% and 20\% of
the total variance of $X_{ij}$, respectively.

%

%
\begin{table}
\def\arraystretch{.92}
\tabcolsep=5pt
\caption{Estimation errors for different sample sizes ($\lambda
_1=0.4,\lambda_2=0.2,\alpha_1=1,\alpha_2=-0.5$)}
\label{tabLa1}
\begin{tabular*}{\tablewidth}{@{\extracolsep{\fill}}lcccd{2.2}ccd{2.2}d{2.2}@{}}
\hline
\multicolumn{2}{@{}l}{\textbf{Sample sizes}}
& \multicolumn{3}{c}{\textbf{Parameter
estimates}} & \multicolumn{4}{c@{}}{\textbf{Prediction}}\\[-4pt]
\multicolumn{2}{@{}l}{\hrulefill}
& \multicolumn{3}{c}{\hrulefill}
& \multicolumn{4}{c@{}}{\hrulefill}\\
$\bolds{n}$ & \multicolumn{1}{c}{$\bolds{p}$}
& \multicolumn{1}{c}{$\bolds{\sum|\widehat\alpha_r-\alpha_r|}$}
& \multicolumn{1}{c}{$\bolds{\sum|\widehat\beta
_r-\beta
_r|}$} & \multicolumn{1}{c}{\textbf{Opt.} $\bolds{\rho}$}
& \multicolumn{1}{c}{\textbf{Sample}}
& \multicolumn{1}{c}{\textbf{Exact}}
& \multicolumn{1}{c}{\textbf{Opt.} $\bolds{\rho}$}
& \multicolumn{1}{c@{}}{$\bolds{C_P}$} \\
\hline
\multicolumn{9}{@{}c@{}}{Lasso applied to augmented model:}\\[4pt]
\hphantom{00,}50 & \hphantom{00,}50 & 0.3334 &0.8389 & 4.53 & 0.0498 & 0.1004 & 1.76 & 1.55 \\
\hphantom{0,}100 & \hphantom{0,}100 & 0.2500 & 0.5774 &6.84 & 0.0328 & 0.0480 & 3.50 & 3.29 \\
\hphantom{0,}250 & \hphantom{0,}250 & 0.1602 & 0.3752 & 12.27 & 0.0167 & 0.0199 & 7.55 & 7.22 \\
\hphantom{0,}500 & \hphantom{0,}500 & 0.1150 & 0.2752 & 18.99 & 0.0096 & 0.0106 & 12.66 & 12.21
\\
5,000 & \hphantom{0,}100 & 0.0378 & 0.0733 & 48.33 & 0.0152 & 0.0154 & 27.48 & 26.74
\\
\hphantom{0,}100 & 2,000 & 0.2741 & 0.8664 & 10.58 & 0.0420 & 0.0651 & 5.42 & 5.24 \\
\hline
$\bolds{n}$ & $\bolds{p}$ &
\multicolumn{1}{c}{$\bolds{\sum|\widehat\beta_r-\beta_r^{\mathit{LR}}|}$}
& \multicolumn{1}{c}{$\bolds{\sum|\widehat\beta
_r-\beta_r|}$} & & \multicolumn{1}{c}{\textbf{Sample}}
& \multicolumn{1}{c}{\textbf{Exact}}
& \multicolumn{1}{c}{\textbf{Opt.} $\bolds{\rho}$} & \\
\hline
\multicolumn{9}{@{}c@{}}{Lasso applied to standard linear regression
model:}\\[4pt]
\hphantom{00,}50 & \hphantom{00,}50 & 2.2597 & 1.8403 & & 0.0521 & 0.1370 & 0.92 & \\
\hphantom{0,}100 & \hphantom{0,}100 & 2.2898 & 1.9090 & & 0.0415 & 0.0725 & 1.90 & \\
\hphantom{0,}250 & \hphantom{0,}250 & 2.3653 & 1.7661 & & 0.0257 & 0.0345 & 4.12 & \\
\hphantom{0,}500 & \hphantom{0,}500 & 2.4716 & 1.7104 & & 0.0174 & 0.0207 & 6.87 & \\
5,000 & \hphantom{0,}100 & 0.5376 & 1.5492 & & 0.0161 & 0.0168 & 10.14 & \\
\hphantom{0,}100 & 2,000 & 3.7571 & 2.2954 & & 0.0523 & 0.1038 & 3.17 & \\
\hline
\end{tabular*}
\vspace*{-12pt}
\end{table}

Figure \ref{figur1} shows estimation results of one typical simulation with
$n=p=100$.
The left panel contains the parameter estimates for the augmented
model. The paths
of estimated
coefficients $\widehat\beta_j$ for the 4 significant variables (black
lines), the 96 variables with $\beta_j=0$ (red lines), as well as of
the untransformed estimates
$\widehat{\widetilde{\alpha}}_r$ (blue lines) of $\alpha_r$, $r=1,2$,
are plotted as a function of $\rho$. The four significant coefficients
as well as
$\alpha_1$ and $\alpha_2$ can immediately been identified in the
figure. The
right panel shows a corresponding plot of estimated coefficients when
Lasso is directly applied to the standard regression model (\ref
{linear-model}). As has to be
expected by (\ref{stanmodident}) the necessity of compensating the
effects of
$\alpha_1,\alpha_2$ by a large number of small, nonzero coefficients
generates a
general ``noise level'' which makes it difficult to identify the four
significant
variables in (\ref{simul-2}). The penalties $\rho$ in the figure as
well as
in subsequent tables have to be interpreted in terms of the scaling used
by the LARS-algorithm and have to be multiplied with $2/n$ in order to
correspond to the standardization used in the proceeding sections.

The upper part of Table \ref{tabLa1} provides simulation results with
respect to the augmented model for different sample sizes $n$ and $p$.
In order to access the quality of parameter estimates we evaluate
$\sum_{r=1}^2 |\widehat\alpha_r-\alpha_r|$ as well as $\sum_{j=1}^p
|\widehat\beta_r-\beta_r|$ at the optimal value of $\rho$, where the
minimum of $\sum_{r=1}^2 |\widehat\alpha_r-\alpha_r|+\sum_{j=1}^p
|\widehat\beta_r-\beta_r|$ is obtained. Moreover, we record the value
of $\rho$ where the minimal sample prediction error
%
%
\begin{equation}\label{simul-3}
\frac{1}{n}\sum_{i=1}^n \Biggl( \sum_{r=1}^2\alpha_r
\xi_{ir}+\sum_{j=1}^p \beta_j X_{ij}- \Biggl(
\sum_{r=1}^2\widehat{\alpha}_r \widehat{\xi}_{ir}+\sum_{j=1}^p
\widehat{\beta}_j X_{ij}\Biggr)\Biggr)^2
\end{equation}
ia attained. For the same value $\rho$ we also determine the exact
prediction error
%
%
\begin{equation}\label{simul-4}\quad
\mathbb{E} \Biggl( \sum_{r=1}^2\alpha_r \xi_{n+1,r}+\sum_{j=1}^p
\beta_j X_{n+1,j}- \Biggl( \sum_{r=1}^2\widehat{\alpha}_r
\widehat{\xi}_{n+1,r}+\sum_{j=1}^p \widehat{\beta}_j
X_{n+1,j}\Biggr)\Biggr)^2
\end{equation}
for a new observation $\mathbf{X}_{n+1}$ independent of $X_1,\ldots
,X_n$. The columns
of Table \ref{tabLa1} report the average values of the corresponding quantities
over the 1,000
replications. To get some insight into a practical choice of the
penalty, the last
column additionally yields the average value of the parameters $\rho$ minimizing
the $C_P$-statistics. $C_P$~is computed by using the R-routine
``summary.lars'' and plugging in the true error variance $\sigma
^2=0.1$. We see\vadjust{\goodbreak} that in all situations
the average value of $\rho$ minimizing $C_P$ is very close to the
average $\rho$
providing the smallest prediction error. The penalties for optimal
parameter estimation are, of course, larger.

It is immediately seen that the quality of estimates considerably
increases when going from $n=p=50$ to $n=p=500$. An interesting result
consists of the fact that the prediction error is smaller for $n=p=500$
than for $n=5\mbox{,}000$, $p=100$. This may be interpreted as a consequence of
(\ref{thmiden-4}).

The lower part of Table \ref{tabLa1} provides corresponding simulation results
with respect to Lasso estimates
based on the standard regression model. In addition to the minimal
error $\sum_{r=1}^p
|\widehat\beta_r-\beta_r|$ in estimating\vspace*{1pt} the parameters $\beta_j$ of
(\ref{simul-2})
we present the minimal $L_1$-distance $\sum_{r=1}^p
|\widehat\beta_r-\beta_r^{\mathit{LR}}|$, where $\beta_1^{\mathit{LR}},
\ldots,\beta_p^{\mathit{LR}}$
is the (nonsparse) set
of parameters minimizing the population prediction error. Sample and
exact prediction errors are obtained by straightforward modifications
of (\ref{simul-3}) and (\ref{simul-4}). Quite obviously, no reasonable
parameter estimates are obtained in the cases with $p\geq n$. Only\vspace*{1pt} for
$n=5\mbox{,}000$, $p=100$, the table indicates a comparably small error
$\sum_{r=1}^p
|\widehat\beta_r-\beta_r^{\mathit{LR}}|$. The prediction error shows a somewhat
better behavior. It is, however, always larger than the prediction
error of the augmented model. The relative difference increases with $p$.

It was mentioned in Section \ref{sec4} that a suitable criterion to estimate
the dimension $k$ of an approximate factor model consists in minimizing
(\ref{baicrit}). This criterion proved to work well in our simulation
study. Recall that the true factor dimension is $k=2$. For $n=p=50$ the
average value of the estimate $\widehat k$ determined
by (\ref{baicrit}) is 2.64. In all other situations reported in Table
\ref{tabLa1} an estimate
$\widehat k =2$ is obtained in each of the 1,000 replications.

%
\begin{table}
\tabcolsep=0pt
\def\arraystretch{.92}
\caption{Estimation errors under
different setups ($n=100, p=250$)}
\label{tabLa2}
\begin{tabular*}{\tablewidth}{@{\extracolsep{\fill}}ld{1.2}cd{2.1}ccd{2.2}ccc@{}}
\hline
&&&&\multicolumn{3}{c}{\textbf{Parameter estimates}}
&\multicolumn{3}{c@{}}{\textbf{Prediction}}\\[-4pt]
&&&&\multicolumn{3}{c}{\hrulefill} & \multicolumn{3}{c@{}}{\hrulefill}\\
\multicolumn{1}{@{}l}{$\bolds{\lambda_1}$}
& \multicolumn{1}{c}{$\bolds{\lambda_2}$}
& \multicolumn{1}{c}{$\bolds{\alpha_1}$}
& \multicolumn{1}{c}{$\bolds{\alpha_2}$} &
\multicolumn{1}{c}{$\bolds{\sum|\widehat\alpha_r-\alpha_r|}$}
& \multicolumn{1}{c}{$\bolds{\sum|\widehat\beta_r-\beta_r|}$}
& \multicolumn{1}{c}{\textbf{Opt.} $\bolds{\rho}$} &
\multicolumn{1}{c}{\textbf{Sample}} &
\multicolumn{1}{c}{\textbf{Exact}} &
\multicolumn{1}{c@{}}{\textbf{Opt.} $\bolds{\rho}$} \\
\hline
\multicolumn{10}{@{}c@{}}{Lasso applied to augmented model:}\\[4pt]
0.06 & 0.03 & 1 & -0.5 & 0.3191 & 0.6104 & 10.37 & 0.0670 & 0.2259 &
2.46 \\
0.2 & 0.1 & 1 & -0.5 & 0.2529 & 0.5335 & 8.19 & 0.0414 & 0.0727 & 3.92
\\
0.4 & 0.2 & 1 & -0.5 & 0.2500 & 0.6498 &7.86 & 0.0319 & 0.0454 & 4.35
\\
0.6 & 0.3 & 1 & -0.5 & 0.2866 & 1.1683 &8.46 & 0.0273 & 0.0350 & 4.56
\\
0.06 & 0.03 & 0 & 0 & 0.0908 & 0.4238 & 7.42 & 0.0257 & 0.0311 & 4.62
\\
0.2 & 0.1 & 0 & 0 & 0.1044 & 0.4788 & 7.64 & 0.0257 & 0.0316 & 4.69 \\
0.4 & 0.2 & 0 & 0 & 0.1192 & 0.6400 &7.84 & 0.0250 & 0.0314 & 4.74 \\
0.6 & 0.3 & 0 & 0 & 0.1825 & 1.1745 &8.78 & 0.0221 & 0.0276 & 5.13 \\
\hline
\multicolumn{1}{@{}l}{$\bolds{\lambda_1}$}
& \multicolumn{1}{c}{$\bolds{\lambda_2}$}
& \multicolumn{1}{c}{$\bolds{\alpha_1}$}
& \multicolumn{1}{c}{$\bolds{\alpha_2}$}
& \multicolumn{1}{c}{$\bolds{\sum|\widehat\beta
_r-\beta_r^{\mathit{LR}}|}$}
& \multicolumn{1}{c}{$\bolds{\sum|\widehat\beta_r-\beta_r|}$}
& & \multicolumn{1}{c}{\textbf{Sample}}
& \multicolumn{1}{c}{\textbf{Exact}} &
\multicolumn{1}{c@{}}{\textbf{Opt.} $\bolds{\rho}$} \\
\hline
\multicolumn{10}{@{}c@{}}{Lasso applied to standard linear regression
model:}\\[4pt]
0.06 & 0.03 & 1 & -0.5 & 5.0599 & 1.9673 & & 0.0777 & 0.3758 & 1.95 \\
0.2 & 0.1 & 1 & -0.5 & 3.4465 & 2.3662 & & 0.0583 & 0.1403 & 2.63 \\
0.4 & 0.2 & 1 & -0.5 & 2.9215 & 2.0191 & & 0.0425 & 0.0721 & 2.45 \\
0.6 & 0.3 & 1 & -0.5 & 3.2014 & 2.2246 & & 0.0277 & 0.0387 & 1.47 \\
0.06 & 0.03 & 0 & 0 & 0.4259 & 0.4259 & & 0.0216 & 0.0285 & 4.80 \\
0.2 & 0.1 & 0 & 0 & 0.4955 & 0.4955 & & 0.0222 & 0.0295 & 4.24 \\
0.4 & 0.2 & 0 & 0 & 0.6580 & 0.6580 & & 0.0228 & 0.0303 & 3.17 \\
0.6 & 0.3 & 0 & 0 & 1.1990 & 1.1990 & & 0.0215 & 0.0283 & 1.66 \\
\hline
\end{tabular*}
\vspace*{-12pt}
\end{table}

Finally, Table \ref{tabLa2} contains simulations results for $n=100$, $p=250$, and
different values of
$\lambda_1,\lambda_2,\alpha_1,\alpha_2$. All columns have to be
interpreted similar to those of Table \ref{tabLa1}. For $\alpha_1=1,\alpha_2=-0.5$
suitable parameter estimates can obviously only been determined by
applying the augmented model. For
$\alpha_1=\alpha_2=0$ model (\ref{simul-2}) reduces to a sparse,
standard linear
regression model. It is then clearly unnecessary to apply the
augmented model.
Both methods then lead to roughly equivalent parameter estimates.

We want to emphasize that $\lambda_1=0.06$, $\lambda_2=0.03$
constitutes a particularly
difficult situation. Then the first and second factor only explain
6\% and 3\% of the variance of $X_{ij}$. Consequently, $v(k)$ is very
small and one will expect a fairly large error in estimating $\xi
_{ir}$. Somewhat surprisingly the augmented
model still provides reasonable parameter estimates,
the only problem in this case seems to be a fairly large prediction error.

Another difficult situation in an opposite direction is $\lambda
_1=0.6$, $\lambda_2=0.3$. Then both factors together explain 90\%
of 
the variability of $X_{ij}$, while $Z_{ij}$ only explains the remaining
10\%. Consequently, $D_1$ is very small and one may expect problems in
the context of the restricted eigenvalue condition. The
table shows that this case yields the smallest
prediction error, but the quality of parameter estimates deteriorates.

\begin{appendix}
\section*{Appendix}\label{app}

\begin{pf*}{Proof of Proposition \ref{propcov}}
Define $Q_{ijl}=V_{ij}V_{il}-\mathbb{E}(V_{ij}V_{il})$, $i=1,\ldots,n$,
$1\le j,l\le p$. For any $C>0$ and $\varepsilon>0$, noting that
$\mathbb{E}(Q_{ijl})=0$, we have
\begin{eqnarray*}
&&
P\Biggl( \Biggl|\frac{1}{n}\sum_{i=1}^nV_{ij}V_{il}-\mathbb{E}(V_{ij}V_{il})\Biggr|
> \varepsilon\Biggr) \\
&&\qquad= P\Biggl(\Biggl|\frac{1}{n}\sum_{i=1}^nQ_{ijl}\Biggr| >
\varepsilon\Biggr)\\
&&\qquad= P\Biggl(\Biggl|\frac{1}{n}\sum_{i=1}^nQ_{ijl}I(|Q_{ijl}|\le C)-\mathbb
{E}\bigl(Q_{ijl}I(|Q_{ijl}|\le C)\bigr)\\
&&\qquad\quad\hspace*{30.5pt}{} +Q_{ijl}I(|Q_{ijl}|> C)-\mathbb
{E}\bigl(Q_{ijl}I(|Q_{ijl}|> C)\bigr)\Biggr| > \varepsilon\Biggr)\\
&&\qquad\le P\Biggl(\Biggl|\frac{1}{n}\sum_{i=1}^nQ_{ijl}I(|Q_{ijl}|\le
C)-\mathbb
{E}\bigl(Q_{ijl}I(|Q_{ijl}|\le C)\bigr)\Biggr| > \varepsilon/2 \Biggr)\\
&&\qquad\quad{} + P\Biggl(\Biggl|\frac{1}{n}\sum_{i=1}^nQ_{ijl}I(|Q_{ijl}|> C)-\mathbb
{E}\bigl(Q_{ijl}I(|Q_{ijl}|> C)\bigr)\Biggr| > \varepsilon/2 \Biggr),
\end{eqnarray*}
where $I(\cdot)$ is the indicator function. We have
\[
\bigl|Q_{ijl}I(|Q_{ijl}|\le C)-\mathbb{E}\bigl(Q_{ijl}I(|Q_{ijl}|\le
C)\bigr)\bigr|\le2C
\]
and
\begin{eqnarray*}
\mathbb{E}\bigl(\bigl(Q_{ijl}I(|Q_{ijl}|\le C)-\mathbb{E}\bigl(Q_{ijl}I(|Q_{ijl}|\le
C)\bigr)\bigr)^2\bigr)&\le& \operatorname{Var}(V_{ij}V_{il})\\
&\le&(\mathbb{E}(V_{ij}^4)\mathbb{E}(V_{il}^4))^{1/2}\\
&\le& C_1.
\end{eqnarray*}
Applying the Bernstein inequality for bounded centered random variables
[see \citet{Ho63}] we get
%
%
\setcounter{equation}{0}
\begin{eqnarray}\label{term1}
&&P\Biggl(\Biggl|\frac{1}{n}\sum_{i=1}^nQ_{ijl}I(|Q_{ijl}|\le C)-\mathbb
{E}\bigl(Q_{ijl}I(|Q_{ijl}|\le C)\bigr)\Biggr| > \frac{\varepsilon}{2}
\Biggr)\nonumber\\[-8pt]\\[-8pt]
&&\qquad\le\exp\biggl\{ \frac{-\varepsilon^2n}{8(C_1+C\varepsilon/3)}
\biggr\} .\nonumber
\end{eqnarray}
We have now
%
%
\begin{eqnarray}\label{term2}
&&P\Biggl(\Biggl|\frac{1}{n}\sum_{i=1}^nQ_{ijl}I(|Q_{ijl}|> C)-\mathbb
{E}\bigl(Q_{ijl}I(|Q_{ijl}|> C)\bigr)\Biggr| > \varepsilon/2
\Biggr)\nonumber\\[-8pt]\\[-8pt]
&&\qquad\le\sum_{i=1}^nP(|Q_{ijl}|> C)+P\bigl(\mathbb
{E}\bigl(|Q_{ijl}|I(|Q_{ijl}|> C)\bigr) > \varepsilon/4
\bigr).\nonumber
\end{eqnarray}
Using Markov's inequality and (\ref{exponential-condition}) we obtain
%
%
\begin{eqnarray}\label{term3}
P(|Q_{ijl}|> C)&\le& P\bigl(|V_{ij}|> \sqrt{C/2}\bigr)+P\bigl(|V_{il}|> \sqrt
{C/2}\bigr)\nonumber\\
&&{} +P\bigl((\mathbb{E}(V_{ij}^2)\mathbb
{E}(V_{il}^2))^{1/2}>C/2\bigr)\\
&\le& \frac{2C_1}{e^{a\sqrt{C/2}}}+P\bigl((\mathbb{E}(V_{ij}^2)\mathbb
{E}(V_{il}^2))^{1/2}>C/2\bigr).\nonumber
\end{eqnarray}
Choose $\varepsilon=C_0\sqrt{\log p/n}$ and $C=\sqrt{C_0n/\log p}$, where
$C_0$ is a positive constant such that $C_1^{1/2}\le\frac{1}{2}\sqrt
{C_0 n\log p}$ and $C_1\le\frac{1}{8}C_0e^{a\sqrt{C_0n/\log p}}\sqrt
{\log p/n}$. Note now that
\[
P\bigl((\mathbb{E}(V_{ij}^2)\mathbb{E}(V_{il}^2))^{1/2}>C/2\bigr)=0,
\]
while
\[
\mathbb{E}\bigl(|Q_{ijl}|I(|Q_{ijl}|> C)\bigr)\le(\mathbb{E}(V_{ij}^2)\mathbb
{E}(V_{il}^2))^{1/2}P(|Q_{ijl}|> C)\le\frac{2C_1^{3/2}}{e^{a\sqrt
{C_0n/\log p}}},
\]
which implies
\[
P\bigl(\mathbb{E}\bigl(|Q_{ijl}|I(|Q_{ijl}|> C)\bigr) > \varepsilon/4 \bigr)=0.
\]
Inequalities (\ref{term1}), (\ref{term2}) and (\ref{term3}) lead
finally to
%
%
\begin{eqnarray}\label{term4}
&&P\Biggl( \Biggl|\frac{1}{n}\sum_{i=1}^nV_{ij}V_{il}-\mathbb{E}(V_{ij}V_{il})\Biggr|
> C_0\sqrt{\frac{\log p}{n}} \Biggr)\nonumber\\[-8pt]\\[-8pt]
&&\qquad\le p^{-C^2_0/(8(C_1+C_0^{3/2}/3))}+2nC_1e^{-({a/2})(n/\log
p)^{1/4}}.\nonumber
\end{eqnarray}
The result (\ref{propcovresult}) is now a consequence of (\ref
{term4}) since
\begin{eqnarray*}
&&P\Biggl(\sup_{1\le j,l\le p} \Biggl|\frac{1}{n}\sum
_{i=1}^nV_{ij}V_{il}-\mathbb{E}(V_{ij}V_{il})\Biggr| > C_0\sqrt{\frac{\log
p}{n}} \Biggr)\\
&&\qquad\le\sum_{j=1}^p\sum_{l=1}^pP\Biggl( \Biggl|\frac{1}{n}\sum
_{i=1}^nV_{ij}V_{il}-\mathbb{E}(V_{ij}V_{il})\Biggr| > C_0\sqrt{\frac{\log
p}{n}} \Biggr).
\end{eqnarray*}
\upqed
\end{pf*}
\begin{pf*}{Proof of Theorem \ref{thmiden1}}
For any symmetric matrix $\mathbf
{A}$, we denote by $\lambda_1(\mathbf{A})>\lambda_2(\mathbf{A})>\cdots$
its eigenvalues. Weyl's
perturbation theorem [see, e.g., \citet{bel61}, page 63] implies that
for any symmetric matrices $\mathbf{A}$ and $\mathbf{B}$
and all $r=1,2,\ldots$
%
%
\begin{equation}\label{knut01}
| \lambda_r(\mathbf{A}+\mathbf{B})-\lambda_r(\mathbf{A})|
\leq\|\mathbf{B}\|,
\end{equation}
where $\|\mathbf{B}\|$ is the usual matrix norm defined as
\[
\|\mathbf{B}\| = \sup_{\|\mathbf{u}\|_2=1}( \mathbf{u}^T\mathbf
{B}\mathbf{B}^T\mathbf{u})^{1/2}.
\]
Since $\frac{1}{p}\bolds{\Sigma}=\frac{1}{p}\bolds{\Gamma
}+\frac{1}{p}\bolds{\Psi}$,
(\ref{knut01}) leads to $|\mu_r-\lambda_r|\leq\|\frac{1}{p}\bolds
{\Psi}\|$.
By assumption, $\frac{1}{p}\bolds{\Psi}$ is a diagonal matrix
with diagonal entries $\frac{D_1}{p}\leq
\frac{\sigma_i^2}{p}\leq\frac{D_2}{p}$, $j=1,\ldots,p$. Therefore $\|
\frac{1}{p}\bolds{\Psi}\|\leq\frac{D_2}{p}$ and (\ref{thmiden-1})
is an immediate consequence.

In order to verify (\ref{thmiden-2}) first note that Lemma A.1 of
\citet{knUt01} implies that for symmetric matrices $\mathbf{A}$ and
$\mathbf{B}$
%
%
\begin{eqnarray}\label{knut04}
\| \bolds{\psi}_r(\mathbf{A}+\mathbf{B})-\bolds{\psi
}_r(\mathbf{A})\|_2
&\leq&
\frac{\|\mathbf{B}\|}{\min_{j\neq r}|\lambda_j(\mathbf{A})-\lambda
_r(\mathbf{A})|}\nonumber\\[-8pt]\\[-8pt]
&&{}+\frac{6\|\mathbf{B}\|^2}{\min_{j\neq r}|\lambda_j(\mathbf{A})-\lambda
_r(\mathbf{A})|^2},\nonumber
\end{eqnarray}
where $\bolds{\psi}_1(\mathbf{A})$, $\bolds{\psi}_2(\mathbf
{A}),\ldots$ are the eigenvectors corresponding
to the eigenvalues $\lambda_1(\mathbf{A})>\lambda_2(\mathbf{A})>\cdots
.$ By assumption (A.3)
this implies
\[
\|\mu_r-\psi_r\|_2\leq\frac{\|({1}/{p})\bolds{\Psi}\|}{v(k)}
+\frac{6\|({1}/{p})\bolds{\Psi}\|^2}{v(k)^2}\leq\frac{2D_2}{pv(k)}
\]
for all $r=1,\ldots,k$.
Since $D_0\geq\mathbb{E}( X_{ij}^2)=\sum_{r=1}^p \delta_{rj}^2p\mu_r$,
the second part of (\ref{thmiden-2}) follows from
$D_0\geq\delta_{rj}^2p\mu_r$, $j=1,\ldots,p$.

By (A.3) we necessarily have $ \delta_r\geq\lambda_r\geq v(k)$ for all
$r=1,\ldots,k$.
Consequently, $\sqrt{\mu_r}-\sqrt{\lambda_r}=\frac{\mu_r-\lambda
_r}{\sqrt{\mu_r}+\sqrt{\lambda_r}}
\leq\frac{\mu_r-\lambda_r}{2\sqrt{v(k)}}$.
Furthermore,
note that $\xi_{ir}^*=\xi_{ir}+\frac{\bolds{\delta}^T_r\mathbf
{Z}_i}{\sqrt{p\mu_r}}
+(\frac{\bolds{\delta}_r}{\sqrt{p\mu_r}}-\frac{\bolds{\psi
}_r}{\sqrt{p\lambda_r}})^T
\mathbf{W}_i$. Since $W_i$ and $Z_i$ are uncorrelated, (\ref
{thmiden-1}) and (\ref{thmiden-2})
lead to
\begin{eqnarray*}
&&\mathbb{E}([\xi_{ir}-
\xi_{ir}^*]^2)\\
&&\qquad=
\frac{1}{\mu_r} \bolds\delta^T_r \frac{1}{p}\bolds{\Psi}\bolds\delta_r\\
&&\qquad\quad{}+\biggl(\frac{\bolds{\delta}_r-\bolds{\psi}_r}{\sqrt{\mu_r}}
+\frac{\sqrt{\mu_r}-\sqrt{\lambda_r}}{\sqrt{\mu_r}\sqrt{\lambda
_r}}\bolds{\psi}_r\biggr)^T
\frac{1}{p}\bolds{\Gamma}\biggl(\frac{\bolds{\delta}_r-\bolds
{\psi}_r}{\sqrt{\mu_r}}
+\frac{\sqrt{\mu_r}-\sqrt{\lambda_r}}{\sqrt{\mu_r}\sqrt{\lambda
_r}}\bolds{\psi}_r\biggr)\\
&&\qquad\leq\frac{D_2}{p\mu_r}+2\frac{\lambda_1}{\mu_r}\|\bolds{\delta
}_r-\bolds{\psi}_r\|^2
+2\frac{1}{\mu_r}\bigl(\sqrt{\mu_r}-\sqrt{\lambda_r}\bigr)^2\\
&&\qquad\leq\frac{D_2}{p\mu_r}+\frac{8\lambda_1D_2^2}{\mu_rp^2v(k)^2}+\frac
{D_2^2}{2\mu_rp^2v(k)}.
\end{eqnarray*}
Since $\frac{\bolds{\psi}_r}{\sqrt{p\lambda_r}}^T
\mathbf{X}_i=\xi_{ir}+\frac{\bolds{\psi}_r}{\sqrt{p\lambda_r}}^T
\mathbf{Z}_i$ the second part of (\ref{thmiden-3}) follows from similar
arguments.
Finally, using the Cauchy--Schwarz inequality (\ref{thmiden-4}) is a
straightforward
consequence of (\ref{thmiden-3}).
\end{pf*}
\begin{pf*}{Proof of Theorem \ref{thmapp1}}
With $\mathbf{A}=\bolds{\Gamma}$ and $\mathbf{B}=\widehat
{\bolds{\Sigma}}-\bolds{\Gamma}=\widehat{\bolds{\Sigma
}}-\bolds{\Sigma}+\bolds{\Psi}$, inequality (\ref{knut01})
implies that for all $r\in\{1,\ldots,k\}$
%
%
\begin{equation}\label{knut02}
|\widehat\lambda_r-\lambda_r|\leq\biggl\|\frac{1}{p}\bolds{\Psi}+
\frac{1}{p}(\widehat{\bolds{\Sigma}}-\bolds{\Sigma})\biggr\|
\leq\biggl\|\frac{1}{p}\bolds{\Psi}\biggr\|+\biggl\|\frac{1}{p}(\widehat{\bolds
{\Sigma}}-\bolds{\Sigma})\biggr\|.
\end{equation}
But under events (\ref{B21})--(\ref{B24}) we have
\begin{eqnarray*}
&&\biggl\| \frac{1}{p}(\widehat{\bolds{\Sigma}} -
\bolds{\Sigma}) \biggr\| \\
&&\qquad= \sup_{\| \mathbf{u}\|
_2=1}\biggl[ \mathbf{u}^T \frac{1}{p^2}(\widehat{\bolds{\Sigma}}
-\bolds{\Sigma})^2\mathbf{u}\biggr]^{1/2} \\
&&\qquad= \sup_{\| \mathbf{u}\|_2=1}\Biggl[ \frac{1}{p^2}\sum_{j=1}^p
\Biggl(\sum_{l=1}^p\Biggl(\frac{1}{n}\sum
_{i=1}^nX_{i,j}X_{i,l}-\operatorname{Cov}(X_{i,j},X_{i,l})\Biggr)u_l\Biggr)^2 \Biggr]^{1/2}
\\
&&\qquad\leq\sup_{\| \mathbf{u}\|_2=1}\Biggl[ \frac{1}{p^2}\sum_{j=1}^p \|
\mathbf{u}\|_2^2 \sum_{l=1}^p\Biggl(\frac{1}{n}\sum
_{i=1}^nX_{i,j}X_{i,l}-\operatorname{Cov}(X_{i,j},X_{i,l})\Biggr)^2 \Biggr]^{1/2} \\
&&\qquad\leq
C_0 \sqrt{\frac{\log p}{n}}.
\end{eqnarray*}
On the other hand, $\|\frac{1}{p}\bolds{\Psi}\|\le\frac{D_2}{p}$,
and we can conclude that
%
%
\begin{equation}\label{knut03}
\biggl\|\frac{1}{p}\bolds{\Psi}\biggr\|+\biggl\|
\frac{1}{p}(\widehat{\bolds{\Sigma}}-\bolds{\Sigma})\biggr\|
\leq\frac{D_2}{p}+ C_0 \sqrt{\frac{\log p}{n}}.
\end{equation}
Relation (\ref{thmapp1-1}) now is an immediate consequence of (\ref
{knut02}) and (\ref{knut03}).

Relations (\ref{knut04}) and (\ref{knut03})
together with (A.3),
(A.4) and (\ref{B21})--(\ref{B24}) lead to
\begin{eqnarray*}
\|\widehat{\bolds{\psi}}_r-\bolds{\psi}_r\|_2
&\leq&\frac{\|(1/p)\bolds{\Psi}+
(1/p)(\widehat{\bolds{\Sigma}}-\bolds{\Sigma})\|}{\min
_{j\neq l}|\lambda_j-\lambda_l|}\\
&&{} +
\frac{6\|(1/p)\bolds{\Psi}+
(1/p)(\widehat{\bolds{\Sigma}}-\bolds{\Sigma})\|
_2^2}{\min_{j\neq l}|\lambda_j-\lambda_l|^2}\\
&\leq&\frac{{D_2}/{p}+C_0(\log p/n)^{1/2}}{v(k)} \\
&&{} +
\frac{6({D_2}/{p}+C_0(\log p/n)^{1/2})^2}{v(k)^2}\\
&\leq&2\frac{{D_2}/{p}+C_0(\log p/n)^{1/2}}{v(k)},
\end{eqnarray*}
which gives (\ref{thmapp1-2}). It remains to show (\ref{thmapp1-3}) and
(\ref{thmapp1-4}). Note that the spectral decompositions of $\bolds
{\Psi}$ and
$\widehat{\bolds{\Sigma}}$ imply that for all $j=1,\ldots,p$
\[
\mathbb{E}(W_{ij}^2)=\sum_{r=1}^p \psi_{rj}^2p\lambda_r,\qquad
\frac{1}{n}\sum_{i=1}^n X_{ij}^2=\sum_{r=1}^p \widehat{\psi
}_{rj}^2p\widehat{\lambda}_r.
\]
Under events (\ref{B21})--(\ref{B24})
we therefore obtain for all $r\le k$
%
%
\begin{eqnarray}
\label{knut05}
\psi_{rj}^2&\leq&\frac{\mathbb{E}(W_{ij}^2)}{p\lambda_r}\leq\frac
{D_0-D_1}{p\lambda_r}\leq\frac{D_0-D_1}{pv(k)},
\\
\label{knut06}
\widehat{\psi}_{rj}^2&\leq&\frac{({1}/{n})\sum_{i=1}^n
X_{ij}^2}{p\widehat{\lambda}_r}\leq\frac{D_0+C_0(\log
p/n)^{1/2}}{p\widehat{\lambda}_r}.
\end{eqnarray}
But by assumptions (A.3) and (A.4), relation (\ref{thmapp1-1}) leads to
$\widehat{\lambda}_r\geq\frac{5v(k)}{6}$. Equations (\ref{thmapp1-3}) and (\ref
{thmapp1-4})
then are immediate consequences of (\ref{knut05}) and (\ref{knut06})
\end{pf*}
\begin{pf*}{Proof of Theorem \ref{thmapp2}}
Choose an arbitrary $j\in\{1,\ldots
,p\}$. Note that
$(\widehat{\mathbf{P}}_k\mathbf{X}_i)_j=X_{ij}-\sum_{r=1}^k \widehat
{\psi}_{rj}\widehat{\bolds{\psi}}{}^T_r \mathbf{X}_i$.
Since $X_{ij}=W_{ij}+Z_{ij}$
we obtain the decomposition
%
%
\begin{eqnarray}\label{prthmapp10}
\frac{1}{n}\sum_{i=1}^n (\widehat{\mathbf{P}}_k\mathbf{X}_i)_j^2&=&
\frac{1}{n}\sum_{i=1}^n \Biggl(W_{ij}-\sum_{r=1}^k \widehat{\psi}_{rj}\widehat
{\bolds{\psi}}^T_r \mathbf{X}_i\Biggr)^2 \nonumber\\
&&{}+2 \frac{1}{n}\sum_{i=1}^n Z_{ij}\Biggl(W_{ij}-\sum_{r=1}^k \widehat{\psi
}_{rj}\widehat{\bolds{\psi}}{}^T_r \mathbf{X}_i\Biggr)\\
&&{}
+\frac{1}{n}\sum_{i=1}^n Z_{ij}^2.
\nonumber
\end{eqnarray}
Under
events (\ref{B21})--(\ref{B24}), we have $|\sigma_j^2-\frac{1}{n}\sum
_{i=1}^n Z_{ij}^2|\leq C_0n^{-1/2}\sqrt{\log p}$ as well as
$|\frac{1}{n}\sum_{i=1}^n Z_{ij}W_{ij}|\le C_0n^{-1/2}\sqrt{\log p}$.
Furthermore,
$\mathbb{E}(Z_{ij}X_{ij})=\sigma^2_j$ and $\mathbb{E}(Z_{ij}X_{il})=0$
for $j\neq l$. Therefore,
\begin{eqnarray*}
\Biggl|\sum_{r=1}^k \widehat{\psi}_{rj}\Biggl(\frac{1}{n}\sum_{i=1}^n Z_{ij}\widehat
{\bolds{\psi}}^T_r \mathbf{X}_i\Biggr)\Biggr|
\leq\sum_{r=1}^k \widehat{\psi}_{rj}^2\sigma_j^2 +
2C_0\sqrt{\frac{\log p}{n}}\sum_{r=1}^k |\widehat{\psi}_{rj}| \Biggl(\sum
_{l=1}^p |\widehat{\psi}_{rl}|\Biggr).
\end{eqnarray*}
Obviously, $\sum_{l=1}^p |\widehat{\psi}_{rl}|\leq\sqrt{p\sum_{l=1}^p
\widehat{\psi}_{rl}^2}=\sqrt{p}$. It
now follows from Theorem \ref{thmapp1}
that there exists a constant $M_1<\infty$, which can be chosen
independently of
all values $n,p,k,S$
satisfying assumptions (A.3) and (A.4), such that
%
%
\begin{eqnarray}\label{prthmapp11}
\frac{1}{n}\sum_{i=1}^n (\widehat{\mathbf{P}}_k\mathbf{X}_i)_j^2-\sigma
_j^2&\geq&
\frac{1}{n}\sum_{i=1}^n \Biggl(W_{ij}-\sum_{r=1}^k
\widehat{\psi}_{rj}\widehat
{\bolds{\psi}}^T_r \mathbf{X}_i\Biggr)^2\nonumber\\
&&{}
- M_1 \frac{k}{v(k)^{1/2}}\sqrt{\frac{\log p}{n}} \\
&\geq& -M_1 \frac{k}{v(k)^{1/2}}\sqrt{\frac{\log p}{n}}.
\nonumber
\end{eqnarray}
Since events (\ref{B21})--(\ref{B24}) have probability $A(n,p)$,
assertion (\ref{pkapp0}) is an immediate consequence.

In order to show (\ref{pkapp1}) first recall that the eigenvectors of
$\frac{1}{p}\widehat{\bolds{\Sigma}}$ possess the well-known
``best basis'' property, that is,
\begin{eqnarray*}
\frac{1}{n}\sum_{i=1}^n\|\widehat{\mathbf{P}}_k\mathbf{X}_i\|_2^2&=&
\frac{1}{n}\sum_{i=1}^n
\Biggl\|\mathbf{X}_i-\sum_{r=1}^k \widehat{\bolds{\psi}}_{r}(\widehat
{\bolds{\psi}}^T_r \mathbf{X}_i)\Biggr\|_2^2\\
&=&
\min_{\mathbf{w}_1,\ldots,\mathbf{w}_k\in\mathbb{R}^p}
\frac{1}{n}\sum
_{i=1}^n \min_{\theta_1,\ldots,\theta_k\in\mathbb{R}}
\Biggl\|X_i-\sum_{r=1}^k \theta_r \mathbf{w}_r\Biggr\|_2^2.
\end{eqnarray*}
For $j=1,\ldots,p$ and $r=1,\ldots,k$
define $\widetilde{\bolds{\psi}}_r^{(j)}\in\mathbb{R}^p$ by
$\widetilde{\psi}_{rj}^{(j)}=\psi_{rj}$ and
$\widetilde{\psi}_{rl}^{(j)}=\widehat{\psi}_{rl}$, $l\neq j$. The
above property then implies
that for any $j$
\[
\frac{1}{n}\sum_{i=1}^n\|\widehat{\mathbf{P}}_k\mathbf{X}_i\|_2^2\leq
\frac{1}{n}\sum_{i=1}^n
\Biggl\|X_i-\sum_{r=1}^k \widetilde{\bolds{\psi}}_{r}^{(j)}(\widehat
{\bolds{\psi}}^T_r \mathbf{X}_i)\Biggr\|_2^2.
\]
Since the vectors $\widehat{\mathbf{P}}_k\mathbf{X}_i$ and
$\mathbf{X}_i-\sum_{r=1}^p \widetilde{\bolds{\psi
}}_{r}^{(j)}(\widehat{\bolds{\psi}}{}^T_r \mathbf{X}_i)$ only
differ in the $j$th element, one can conclude that for any $j=1,\ldots,p$
%
%
\begin{equation}\label{prthmapp1}
\frac{1}{n}\sum_{i=1}^n (\widehat{\mathbf{P}}_k\mathbf{X}_i)_j^2
\leq\frac{1}{n}\sum_{i=1}^n
\Biggl(X_{ij}-\sum_{r=1}^k \psi_{rj} (\widehat{\bolds{\psi}}{}^T_r \mathbf
{X}_i)\Biggr)^2.
\end{equation}
The spectral decomposition of $\widehat{\bolds{\Sigma}}$ implies that
$\widehat{\bolds{\Sigma}}=\sum_{r=1}^p p\widehat{\lambda}_r
\widehat{\bolds{\psi}}_r\widehat{\bolds{\psi}}{}^T_r$ with
\[
p\widehat{\lambda}_r=\frac{1}{n}\sum_{i=1}^n (\widehat{\bolds{\psi
}}^T_r \mathbf{X}_i)^2,\qquad
\frac{1}{n}\sum_{i=1}^n (\widehat{\bolds{\psi}}{}^T_r \mathbf
{X}_i)(\widehat{\bolds{\psi}}_s^T \mathbf{X}_i)=0,\qquad s\neq r.
\]
It therefore follows from (\ref{prthmapp1}) that
%
%
\begin{eqnarray}\label{prthmapp2}
\frac{1}{n}\sum_{i=1}^n (\widehat{\mathbf{P}}_k\mathbf{X}_i)_j^2
&\leq&\frac{1}{n}\sum_{i=1}^n X_{ij}^2 -2\sum_{r=1}^k \psi_{rj} \widehat
{\bolds{\psi}}^T_r \Biggl(\frac{1}{n}\sum_{i=1}^n X_{ij} \mathbf{X}_i\Biggr)
\nonumber\\[-8pt]\\[-8pt]
&&{}+\sum_{r=1}^k \psi_{rj}^2 p\widehat{\lambda}_r.\nonumber
\end{eqnarray}
We obtain $\mathbb{E}(X_{ij}^2)=\mathbb{E}(W_{ij}^2)+\sigma_j^2$ as
well as $\mathbb{E}(X_{ij}X_{il})=\mathbb{E}(W_{ij}W_{il})$ for $j\neq l$.
At the same time under
events (\ref{B21})--(\ref{B24}),
\begin{eqnarray*}
&&\Biggl|\sum_{r=1}^k \psi_{rj} \widehat{\bolds{\psi}}{}^T_r \Biggl(\frac
{1}{n}\sum_{i=1}^n X_{ij} \mathbf{X}_i\Biggr)
-\sum_{r=1}^k \psi_{rj} \bolds{\psi}^T_r \mathbb{E}( W_{ij}
\mathbf{W}_i)\Biggr|\\
&&\qquad\leq C_0n^{-1/2}\sqrt{\log p}\sum_{r=1}^k |\psi_{rj}| \Biggl(\sum_{l=1}^p
|\widehat{\psi}_{rl}|\Biggr)\\
&&\qquad\quad{}+\sum_{r=1}^k |\psi_{rj}||(\bolds{\psi}_r-\widehat{\bolds
{\psi}}_r)^T\mathbb{E}( W_{ij} \mathbf{W}_i)|\\
&&\qquad\quad{}+
\sum_{r=1}^k |\psi_{rj}||\widehat{\psi}_{rj}|\sigma_j^2.
\end{eqnarray*}
Note that $\mathbb{E}( W_{ij}W_{ik})\leq D_0-D_1$ for all $j,k=1,\ldots,p$.
By the Cauchy--Schwarz inequality, Theorem \ref{thmapp1} and
Assumption (A.4), we have
\begin{eqnarray*}
|(\bolds{\psi}_r-\widehat{\bolds{\psi}}_r)^T\mathbb{E}(
W_{ij} \mathbf{W}_i)|
&\leq&\|\bolds{\psi}_r-\widehat{\bolds{\psi}}_r\|_2\|\mathbb
{E}( W_{ij} \mathbf{W}_i)\|_2\\
&\leq&\frac{10C_0n^{-1/2}\sqrt{\log p}}{v(k)}\sqrt{p(D_0-D_1)}
\end{eqnarray*}
as well as
$\sum_{l=1}^p |\widehat{\psi}_{rl}|\leq\sqrt{p\sum_{l=1}^p \widehat
{\psi}_{rl}^2}=\sqrt{p}$.
The bounds for $\psi_{rj}$ and $\widehat{\psi}_{rl}$ derived in
Theorem \ref{thmapp1} then imply that under
events (\ref{B21})--(\ref{B24}) there exists a constant $\widetilde
{M}_2<\infty$, which can be chosen independently of
all values $n,p,k,S$
satisfying assumptions (A.3) and (A.4),
such that
%
%
\begin{eqnarray}\label{prthmapp2b}
&&\Biggl|\sum_{r=1}^k \psi_{rj} \widehat{\bolds{\psi}}{}^T_r \Biggl(\frac
{1}{n}\sum_{i=1}^n X_{ij} \mathbf{X}_i\Biggr)
-\sum_{r=1}^k \psi_{rj} \bolds{\psi}^T_r \mathbb{E}( W_{ij}
\mathbf{W}_i)\Biggr|\nonumber\\[-8pt]\\[-8pt]
&&\qquad\leq\widetilde{M}_2 \frac{k}{v(k)^{3/2}} n^{-1/2}\sqrt{\log
p}.\nonumber
\end{eqnarray}
At the same time,
by Theorem \ref{thmapp1} it follows that there exist constants
${\widetilde{M}}^{*}_2,\widetilde{M}_2^{**}<\infty$ such that
and $r=1,\ldots,k$,
%
%
\begin{eqnarray}\label{prthmapp3}
|p\widehat{\lambda}_r-p\lambda_r|&\leq&\widetilde{M}_2^* pn^{-1/2}\sqrt
{\log p},\nonumber\\[-8pt]\\[-8pt]
|\psi_{rj}^2 p\widehat{\lambda}_r-\psi_{rj}^2p\lambda_r|&\leq&\frac{
\widetilde{M}_2^{**}}{v(k)} n^{-1/2}\sqrt{\log p}.\nonumber
\end{eqnarray}
Note\vspace*{1pt} that
$\mathbb{E}((\bolds{\psi}^T_r\mathbf{W}_i)^2)=p\lambda _r$. Under
events (\ref{B21})--(\ref{B24}), we can now conclude from
(\ref{prthmapp2})--(\ref{prthmapp3}) that there exists a constant $
\widetilde{M}_2^{***}<\infty$, which can be chosen independently of all
values $n,p,k,S$ satisfying Assumptions (A.3) and (A.4), such that
%
%
\begin{eqnarray} \label{prthmapp5}
&&\frac{1}{n}\sum_{i=1}^n (\widehat{\mathbf{P}}_k\mathbf{X}_i)_j^2\nonumber\\
&&\qquad\leq
\Biggl|\sigma_j^2 +\mathbb{E}(W_{ij}^2)-2\sum_{r=1}^k \psi_{rj} \bolds
{\psi}^T_r \mathbb{E}( W_{ij} \mathbf{W}_i)+\sum_{r=1}^k \psi_{rj}^2
p\lambda_r\Biggr|\nonumber\\[-8pt]\\[-8pt]
&&\qquad\quad{} +\widetilde{M}_2^{***} \frac{k}{v(k)} n^{-1/2}\bigl(\sqrt{\log
p}\bigr)^2\nonumber\\
&&\qquad=\sigma_j^2+\mathbb{E}((\mathbf{P}_k\mathbf{W}_i)_j^2
)+\widetilde{M}_2^{***} \frac{k}{v(k)} n^{-1/2}\sqrt{\log p}.
\nonumber
\end{eqnarray}
Relations (\ref{prthmapp11}) and (\ref{prthmapp5}) imply that under
(\ref{B21})--(\ref{B24})
%
%
\begin{eqnarray}\label{prthmapp12}
&&\frac{1}{n}\sum_{i=1}^n \Biggl(W_{ij}-\sum_{r=1}^k \widehat{\psi}_{rj}\widehat
{\bolds{\psi}}^T_r \mathbf{X}_i\Biggr)^2\nonumber\\[-8pt]\\[-8pt]
&&\qquad\leq\mathbb{E}((\mathbf{P}_k\mathbf{W}_i)_j^2)+ M_2^* \frac
{k}{v(k)^{3/2}} n^{-1/2}\sqrt{\log p}\nonumber
\end{eqnarray}
holds with $M_2^*\leq M_1+\widetilde{M}_2^{***}$. Since events (\ref
{B21})--(\ref{B24}) have probability $A(n,p)$,
assertion (\ref{pkapp1}) of Theorem \ref{thmapp2} now is an immediate
consequence of (\ref{prthmapp11}), (\ref{prthmapp5})
and (\ref{prthmapp12}).

It remains to show (\ref{pkapp2}).
We have
%
%
\begin{eqnarray} \label{thmaug2-eq3}\quad
&&
\frac{1}{n}\sum_{i=1}^n\Biggl(\sum_{r=1}^k(\widehat{\xi}_{ir}-\xi
_{ir})\alpha_r\Biggr)^2 \nonumber\\
&&\qquad\le \alpha^2_{\mathrm{sum}} \frac{1}{n}\sum_{i=1}^n\sum
_{r=1}^k(\widehat{\xi}_{ir}-\xi_{ir})^2\\
&&\qquad\leq \alpha^2_{\mathrm{sum}}\Biggl( 2\sum_{r=1}^k
\frac{1}{n}\sum_{i=1}^n \Biggl( \frac{\widehat{\bolds{\psi
}}^T_r\mathbf{X}_i}
{\sqrt{p\widehat{\lambda}_r}}-\frac{\bolds{\psi}^T_r\mathbf{X}_i}
{\sqrt{p\lambda_r}}\Biggr)^2+ 2\sum_{r=1}^k \frac{1}{n}\sum_{i=1}^n
\Biggl(\frac{\bolds{\psi}^T_r\mathbf{Z}_i}
{\sqrt{p\lambda_r}}\Biggr)^2\Biggr).\nonumber
\end{eqnarray}
But for all $r=1,\ldots,k$
%
%
\begin{eqnarray}\label{thmaug2-eq31}
\frac{1}{n}\sum_{i=1}^n \biggl( \frac{\widehat{\bolds{\psi
}}^T_r\mathbf{X}_i}
{\sqrt{p\widehat{\lambda}_r}}-\frac{\bolds{\psi}^T_r\mathbf{X}_i}
{\sqrt{p\lambda_r}}\biggr)^2
&\leq&\frac{2}{n}\sum_{i=1}^n \frac{(\sqrt{\widehat{\lambda}_r}-\sqrt
{\lambda_r})^2(\widehat{\bolds{\psi}}{}^T_r
\mathbf{X}_i)^2}{p\lambda_r\widehat{\lambda}_r}\nonumber\\[-8pt]\\[-8pt]
&&{}+\frac{2}{n}\sum_{i=1}^n
\frac{((\bolds{\psi}_r-\widehat{\bolds{\psi}}_r)^T\mathbf
{X}_i)^2}{p\lambda_r},\nonumber
\end{eqnarray}
and Theorem \ref{thmapp1} and assumptions (A.1)--(A.4) imply that under
events (\ref{B21})--(\ref{B24}) there exist some constants
$M_3^*,M_3^{**}<\infty$, which can be chosen independently of
all values $n,p,k,S$
satisfying assumptions (A.3) and (A.4), such that
%
%
\begin{equation} \label{lem2-30}\quad
\frac{1}{n}\sum_{i=1}^n \frac{(\sqrt{\widehat{\lambda}_r}-\sqrt{\lambda
_r})^2(\widehat{\bolds{\psi}}{}^T_r
\mathbf{X}_i)^2}{p\lambda_r\widehat{\lambda}_r}= \frac{(\widehat{\lambda
}_r-\lambda_r)^2}
{(\sqrt{\widehat{\lambda}_r}+\sqrt{\lambda_r})^2\lambda_s}
\le\frac{M_3^*}{v(k)^2}\frac{\log p}{n}
\end{equation}
and
%
%
\begin{equation}\label{lem2-31}\qquad
\frac{1}{n}\sum_{i=1}^n
\frac{((\bolds{\psi}_r-\widehat{\bolds{\psi}}_r)^T\mathbf
{X}_i)^2}{p\lambda_r}\leq\frac{\||\bolds{\psi}_r-\widehat
{\bolds{\psi}}_r\|_2^2}{p\lambda_r}
\frac{1}{n}\sum_{i=1}^n\|\mathbf{X}_i\|_2^2
\le\frac{M_3^{**}}{v(k)^{3}}\frac{\log p}{n}
\end{equation}
hold for all $r=1,\ldots,k$.

Now note that our setup implies that $ \frac{1}{n}\sum_{i=1}^n\mathbb
{E}((\frac{\bolds{\psi}^T_r\mathbf{Z}_i}
{\sqrt{p\lambda_r}})^2)\leq\frac{D_2}{pv(k)}$ and
$\operatorname{Var}(\frac{1}{n}
\sum_{i=1}^n(\frac{\bolds{\psi}^T_r\mathbf{Z}_i}
{\sqrt{p\lambda_r}})^2)\leq\frac{D_3+2D_2^2}{nv(k)^2p^2}$ hold
for all
$r=1,\ldots,k$.
The Chebyshev inequality thus implies
that the event
%
%
\begin{equation}\label{B25}
\frac{1}{n}\sum_{i=1}^n\biggl(\frac{\bolds{\psi}^T_r\mathbf{Z}_i}
{\sqrt{p\lambda_r}}\biggr)^2\leq\frac{D_2+\sqrt{D_3+2D_2^2}}{pv(k)}\qquad
\mbox{for all } r=1,\ldots,k
\end{equation}
holds with probability at least $1-\frac{k}{n}$. We can thus infer from
(\ref{thmapp2-eq3})--(\ref{lem2-31})
that there exists some positive constant $M_3<\infty$, which can be
chosen independently of the values $n,p,k,S$ satisfying (A.3)--(A.5),
such that under events (\ref{B21})--(\ref{B24}) and (\ref{B25})
\begin{eqnarray*}
&&\frac{1}{n}\sum_{i=1}^n\Biggl(\sum_{r=1}^k(\widehat{\xi}_{ir}-\xi
_{ir})\alpha_r\Biggr)^2\\
&&\qquad\leq\alpha^2_{\mathrm{sum}}\biggl(\frac{4k(M_3^*+M_3^{**})}{v(k)^3}\frac{\log p}{n}+
\frac{2k(D_2+\sqrt{D_3+2D_2^2})}{pv(k)}\biggr).
\end{eqnarray*}
Recall that events (\ref{B21})--(\ref{B24}) and (\ref{thmapp2-eq1})
simultaneously hold with probability at least $A(n,p)-(p+k)^{-A^2/8}$,
while (\ref{B25}) is satisfied with
probability at least $1-\frac{k}{n}$. This proves
assertion (\ref{pkapp2}) with $M_3=4(M_3^*+M_3^{**})+2(D_2+\sqrt{D_3+2D_2^2})$.
\end{pf*}
\begin{pf*}{Proof of Proposition \ref{propsparse2}}
Let $\widehat{\mathbf{Q}}_k$ denote the $p\times p$ diagonal matrix
with diagonal entries $1/\sqrt{\frac{1}{n}\sum_{l=1}^n(\widehat{\mathbf
{P}}_k\mathbf{X}_l)^2_1},\ldots,
1/\sqrt{\frac{1}{n}\sum_{l=1}^n(\widehat{\mathbf{P}}_k\mathbf
{X}_l)^2_p}$ and split the $(k+p)$-dimensional vector $\bolds
{\Delta}$ in two vectors $\bolds{\Delta}_1$ and
$\bolds{\Delta}_2$, where $\bolds{\Delta}_1$ is the
$k$-dimensional vector with the $k$ upper components of
$\bolds{\Delta}$, and $\bolds{\Delta}_2$ is the
$p$-dimensional vector with the $p$ lower components of
$\bolds{\Delta}$. Then
\begin{eqnarray*}
\bolds{\Delta}^T\frac{1}{n}\sum_{i=1}^n \bolds{\Phi}_i\bolds
{\Phi}_i^{T} \bolds{\Delta} &=&
\bolds{\Delta}_1^T\bolds{\Delta}_1+\bolds{\Delta}^T_2\frac
{1}{n}\sum_{i=1}^n
\widehat{\mathbf{Q}}_k\widehat{\mathbf{P}}_k\mathbf{X}_i\mathbf{X}_i^{T}
\widehat{\mathbf{P}}_k\widehat{\mathbf{Q}}_k\bolds{\Delta}_2\\
&\geq& \bolds{\Delta}_1^T\bolds{\Delta}_1+\bolds{\Delta}_2^T
\widehat{\mathbf{Q}}_k\widehat{\mathbf{P}}_k\bolds{\Psi}
\widehat{\mathbf{P}}_k\widehat{\mathbf{Q}}_k \bolds{\Delta}_2\\
&&{}+
\bolds{\Delta}^T_2
\widehat{\mathbf{Q}}_k\widehat{\mathbf{P}}_k(
\widehat{\bolds{\Sigma}}-\bolds{\Sigma})
\widehat{\mathbf{P}}_k\widehat{\mathbf{Q}}_k \bolds{\Delta}_2.
\end{eqnarray*}
The matrix $\bolds{\Psi}$ is a diagonal matrix with entries $\sigma
_1^2,\ldots,
\sigma_p^2$, and $\widehat{\bolds{\psi}}{}^T_r \bolds{\Psi
}\widehat{\bolds{\psi}}_s
\leq D_2$ for all $r,s$. Together with the bounds for $\widehat{\psi
}_{rj}$ derived in Theorem \ref{thmapp1} we can conclude that under
(\ref{B21})--(\ref{B24}) there exists a constant
$M_4^*<\infty$, which can be chosen independently of
all values $n,p,k,S$
satisfying assumptions (A.3) and (A.4), such that
\begin{eqnarray*}
\hspace*{-5pt}&&\bolds{\Delta}^T_2\widehat{\mathbf{Q}}_k\widehat{\mathbf
{P}}_k\bolds{\Psi}
\widehat{\mathbf{P}}_k\widehat{\mathbf{Q}}_k\bolds{\Delta}_2\\
\hspace*{-5pt}&&\qquad=\bolds{\Delta}^T_2\widehat{\mathbf{Q}}_k\bolds{\Psi}
\widehat{\mathbf{Q}}_k\bolds{\Delta}_2
-2\sum_{r=1}^k \bolds{\Delta}^T_2\widehat{\mathbf{Q}}_k \widehat
{\bolds{\psi}}_r\widehat{\bolds{\psi}}{}^T_r
\bolds{\Psi}\widehat{\mathbf{Q}}_k\bolds{\Delta}_2\\
\hspace*{-5pt}&&\qquad\quad{}
+\sum_{r=1}^k\sum_{s=1}^k
\widehat{\mathbf{Q}}_k\bolds{\Delta}^T_2 \widehat{\bolds{\psi
}}_r\widehat{\bolds{\psi}}{}^T_r\bolds{\Psi}
\widehat{\bolds{\psi}}_s\widehat{\bolds{\psi}}_s^T \widehat
{\mathbf{Q}}_k\bolds{\Delta}_2
\\
\hspace*{-5pt}&&\qquad\geq
\biggl(\frac{D_1}{\max_j (({1}/{n})\sum_{l=1}^n
(\widehat{\mathbf{P}}_k\mathbf
{X}_l)^2_j)}
- \frac{2k M_4^*+k^2M_4^*}{pv(k)\min_j (({1}/{n})\sum_{l=1}^n(\widehat
{\mathbf{P}}_k\mathbf{X}_l)^2_j)}
\biggr)
\Vert\bolds{\Delta}_2\Vert_2^2.
\end{eqnarray*}
We have
\[
\max_{j} \frac{1}{n}\sum_{i=1}^n(\widehat{\mathbf{P}}_k\mathbf
{X}_i)^2_j\leq
\max_{j} \frac{1}{n}\sum_{i=1}^nX_{ij}^2 \le D_0 +\max_{j}\Biggl|\frac
{1}{n}\sum_{i=1}^nX_{ij}^2-\mathbb{E}(X_{ij}^2)\Biggr|,
\]
and since $D_1\leq D_0$, this leads under (\ref{B21})--(\ref{B24}) to
\begin{eqnarray*}
&&\bolds{\Delta}^T_1\bolds{\Delta}_1+\bolds{\Delta}_2^T
\widehat{\mathbf{Q}}_k\widehat{\mathbf{P}}_k\bolds{\Psi}
\widehat{\mathbf{P}}_k\widehat{\mathbf{Q}}_k\bolds{\Delta}_2\\
&&\qquad\geq\biggl(
\frac{D_1}{D_0+C_0n^{_1/2}\sqrt{\log p}}
- \frac{2k M_4^*+k^2M_4^*}{pv(k)\min_j (({1}/{n})\sum_{l=1}^n(\widehat
{\mathbf{P}}_k\mathbf{X}_l)^2_j)}
\biggr)
\Vert\bolds{\Delta}\Vert_2^2.
\end{eqnarray*}
On the other hand
\begin{eqnarray*}
&&\bolds{\Delta}_2^T
\widehat{\mathbf{Q}}_{k}\widehat{\mathbf{P}}_k(
\widehat{\bolds{\Sigma}}-\bolds{\Sigma})
\widehat{\mathbf{P}}_k \widehat{\mathbf{Q}}_{k} \bolds{\Delta}_2\\
&&\qquad=
(\bolds{\Delta}_{2,J_{0,k+S}}+\bolds{\Delta
}_{2,J_{0,k+S}^C})^T\widehat{\mathbf{Q}}_{k}\widehat{\mathbf{P}}_k(
\widehat{\bolds{\Sigma}}-\bolds{\Sigma})
\widehat{\mathbf{P}}_k\widehat{\mathbf{Q}}_{k} (\bolds{\Delta
}_{2,J_{0,k+S}}+\bolds{\Delta}_{2,J_{0,k+S}^C})\\
&&\qquad=\bolds{\Delta}_{2,J_{0,k+S}}^T\widehat{\mathbf{Q}}_{k}\widehat
{\mathbf{P}}_k(
\widehat{\bolds{\Sigma}}-\bolds{\Sigma})
\widehat{\mathbf{P}}_k\widehat{\mathbf{Q}}_{k} \bolds{\Delta
}_{2,J_{0,k+S}}\\
&&\qquad\quad{}+\bolds{\Delta}_{2,J_{0,k+S}^C}^T\widehat{\mathbf{Q}}_{k}
\widehat{\mathbf{P}}_k(
\widehat{\bolds{\Sigma}}-\bolds{\Sigma})
\widehat{\mathbf{P}}_k\widehat{\mathbf{Q}}_{k} \bolds{\Delta
}_{2,J_{0,k+S}^C}\\
&&\qquad\quad{} +2\bolds{\Delta}_{2,J_{0,k+S}}^T\widehat{\mathbf{Q}}_{k}\widehat{\mathbf{P}}_k(
\widehat{\bolds{\Sigma}}-\bolds{\Sigma})
\widehat{\mathbf{P}}_k\widehat{\mathbf{Q}}_{k}
\bolds{\Delta}_{2,J_{0,k+S}^C},
\end{eqnarray*}
where $\bolds{\Delta}_{2,J_{0,k+S}}$, respectively, $\bolds{\Delta
}_{2,J_{0,k+S}^C}$, is the $p$-dimensional vector with the last $p$
coordinates of $\bolds{\Delta}_{J_{0,k+S}}$, respectively, $
\bolds{\Delta
}_{J_{0,k+S}^C}$. The Cauchy--Schwarz inequality leads to $\|\bolds
{\Delta}_{J_{0,k+S}}\|_1\leq(2(k+S))^{1/2} \|\bolds{\Delta
}_{J_{0,k+S}}\|_2$. Since
$\|\bolds{\Delta}_{J_{0,k+S}^C}\|_1\leq c_0\|\bolds{\Delta
}_{J_{0,k+S}}\|_1$ we have
\begin{eqnarray*}
&&|\bolds{\Delta}_{2,J_{0,k+S}}^T\widehat{\mathbf{Q}}_{k}\widehat
{\mathbf{P}}_k (\widehat{\bolds{\Sigma}}-\bolds{\Sigma
})\widehat{\mathbf{P}}_k\widehat{\mathbf{Q}}_{k} \bolds{\Delta
}_{2,J_{0,k+S}^C}|\\
&&\qquad\leq\max_{j,l} \bigl|\bigl(\widehat{\mathbf{Q}}_{k}\widehat{\mathbf{P}}_k(
\widehat{\bolds{\Sigma}}-\bolds{\Sigma})
\widehat{\mathbf{P}}_k\widehat{\mathbf{Q}}_{k}\bigr)_{j,l}\bigr|\|\bolds
{\Delta}_{2,J_{0,k+S}}\|_1\|\bolds{\Delta}_{2,J_{0,k+S}^C}\|
_1\nonumber\\
&&\qquad\leq\max_{j,l} \bigl|\bigl(\widehat{\mathbf{Q}}_{k}\widehat{\mathbf{P}}_k(
\widehat{\bolds{\Sigma}}-\bolds{\Sigma})\widehat{\mathbf{P}}_k
\widehat{\mathbf{Q}}_{k}\bigr)_{j,l}\bigr|\|\bolds{\Delta}_{J_{0,k+S}}\|_1\|
\bolds{\Delta}_{J_{0,k+S}^C}\|_1\nonumber\\
&&\qquad\leq c_0\max_{j,l} \bigl|\bigl(\widehat{\mathbf{Q}}_{k}\widehat{\mathbf{P}}_k(
\widehat{\bolds{\Sigma}}-\bolds{\Sigma})
\widehat{\mathbf{P}}_k \widehat{\mathbf{Q}}_{k}\bigr)_{j,l}\bigr|\|\bolds
{\Delta}_{J_{0,k+S}}\|_1^2\nonumber\\
&&\qquad\leq
2(k+S)c_0\max_{j,l} \bigl|\bigl(\widehat{\mathbf{Q}}_{k}\widehat{\mathbf{P}}_k(
\widehat{\bolds{\Sigma}}-\bolds{\Sigma})\widehat{\mathbf{P}}_k
\widehat{\mathbf{Q}}_{k}\bigr)_{j,l}\bigr|\|\bolds{\Delta}_{J_{0,k+S}}\|_2^2,
\end{eqnarray*}
and the same upper bound holds for the terms $\bolds{\Delta
}_{2,J_{0,k+S}}^T\widehat{\mathbf{Q}}_{k}\widehat{\mathbf{P}}_k(
\widehat{\bolds{\Sigma}}-\bolds{\Sigma})\widehat{\mathbf{P}}_k
\widehat{\mathbf{Q}}_{k} \times\break\bolds{\Delta}_{2,J_{0,k+S}}$ and $\bolds{\Delta
}_{2,J_{0,k+S}^C}^T\widehat{\mathbf{Q}}_{k}\widehat{\mathbf{P}}_k(
\widehat{\bolds{\Sigma}}-\bolds{\Sigma})\widehat{\mathbf{P}}_k
\widehat{\mathbf{Q}}_{k} \bolds{\Delta}_{2,J_{0,k+S}^C}$ so that
\begin{eqnarray*}
&&\bolds{\Delta}^T_2
\widehat{\mathbf{Q}}_k\widehat{\mathbf{P}}_k(
\widehat{\bolds{\Sigma}}-\bolds{\Sigma})
\widehat{\mathbf{P}}_k\widehat{\mathbf{Q}}_k \bolds{\Delta}_2\\
&&\qquad\leq8(k+S) c_0 \max_{j,l} \bigl|\bigl(\widehat{\mathbf{Q}}_k\widehat
{\mathbf{P}}_k(
\widehat{\bolds{\Sigma}}-\bolds{\Sigma})
\widehat{\mathbf{P}}_k\widehat{\mathbf{Q}}_k\bigr)_{j,l}\bigr|
\|\bolds{\Delta}_{J_{0,k+S}}\|_2^2.
\end{eqnarray*}
Obviously, $
\widehat{\bolds{\psi}}{}^T_r (\widehat{\bolds{\Sigma
}}-\bolds{\Sigma})\widehat{\bolds{\psi}}_s
\leq p \max_{j,l}|\frac{1}{n}\sum
_{i=1}^nX_{i,j}X_{i,l}-\operatorname{Cov}(X_{i,j},X_{i,l}|$ for all $r,s$.
Using Theorem \ref{thmapp1}, one can infer that under
(\ref{B21})--(\ref{B24}),
\begin{eqnarray*}
&&\max_{j,l} \bigl|\bigl(\widehat{\mathbf{Q}}_k\widehat{\mathbf{P}}_k(
\widehat{\bolds{\Sigma}}-\bolds{\Sigma})
\widehat{\mathbf{P}}_k\widehat{\mathbf{Q}}_k\bigr)_{j,l}\bigr|\\
&&\qquad\leq
\frac{\max_{j,l}|({1}/{n})\sum_{i=1}^nX_{i,j}X_{i,l}-\operatorname{Cov}(X_{i,j},X_{i,l})|}
{\min_j (({1}/{n})\sum_{l=1}^n(\widehat{\mathbf{P}}_k\mathbf{X}_l)^2_j)}
\\
&&\qquad\quad{} +
2\sum_{r=1}^k\max_{j,l}\bigl| \bigl(\widehat{\bolds{\psi}}_r\widehat
{\bolds{\psi}}^T_r
(\widehat{\bolds{\Sigma}}-\bolds{\Sigma})\bigr)_{j,l}\bigr|\\
&&\qquad\quad{}+\sum_{r=1}^k\sum_{s=1}^k\max_{j,l}\bigl| \bigl( \widehat{\bolds{\psi
}}_r\widehat{\bolds{\psi}}{}^T_r
(\widehat{\bolds{\Sigma}}-\bolds{\Sigma})
\widehat{\bolds{\psi}}_s\widehat{\bolds{\psi}}_s^T
\bigr)_{j,l}\bigr|\\
&&\qquad\leq\frac{M_4^{**} ({k^2}/{v(k)}) n^{-1/2}\sqrt{\log p}}
{\min_j (({1}/{n})\sum_{l=1}^n(\widehat{\mathbf{P}}_k\mathbf{X}_l)^2_j)},
\end{eqnarray*}
where the constant $M_4^{**}<\infty$ can be chosen independently of
all values $n,p,k,S$
satisfying assumptions (A.3) and (A.4). When combining the above inequalities,
the desired result follows from (A.5) and the bound on\break $\min_{j}
\frac{1}{n}\sum_{i=1}^n(\widehat{\mathbf{P}}_k\mathbf{X}_i)^2_j$ to be obtained
from (\ref{pkapp0})
\end{pf*}
\begin{pf*}{Proof of Theorem \ref{thmaug}}
The first step of the proof consists of showing that under events
(\ref{B21})--(\ref{B24}) the following inequality holds with probability
at least $1-(p+k)^{1-A^2/8}$
%
%
\begin{equation}\label{thmapp2-eq1}
2\biggl\|\frac{1}{n}\bolds{\Phi}^T(\mathbf{Y}
-\bolds{\Phi}\bolds{\theta}) \biggr\|_\infty\le\rho,
\end{equation}
where $\rho=A\sigma\sqrt{\frac{\log(k+p)}{n}}+\frac{M_5\bolds
{\alpha}_{\mathrm{sum}}}{v(k)^{3/2}}\sqrt{\frac{\log p}{n}}$, $A>2\sqrt{2}$ and
$M_5$ is a sufficiently large positive constant.

Since $W_{ij}, Z_{ij}$ and, hence, $\widehat{\xi}_{ir}$ and $\widetilde
{X}_{ij}$ are
independent of the i.i.d. error terms
$\varepsilon_i\sim\mathcal{N}(0,\sigma^2)$, it follows from standard
arguments that
%
%
\begin{equation}\label{thmapp2-eq2}
\sup_{1\le r\le k,1\le j\le p}\Biggl\{ \frac{2}{n}\Biggl|\sum
_{i=1}^n\widehat{\xi}_{ir}\varepsilon_i \Biggr|,\frac{2}{n}\Biggl|\sum
_{i=1}^n\widetilde{X}_{ij}\varepsilon_i \Biggr| \Biggr\} \le A\sigma\sqrt
{\frac{\log(k+p)}{n}}
\end{equation}
holds with probability at least $1-(p+k)^{1-A^2/8}$. Therefore, in order
to prove
(\ref{thmapp2-eq1}) it only remains to show that under events (\ref
{B21})--(\ref{B24})
there exists a positive constant $M_5<\infty$, which can be chosen
independently of the values $n,p,k,S$ satisfying (A.3)--(A.5), such that
%
%
\begin{equation}\label{thmapp2-eq3}
\sup_{1\le r\le k,1\le j\le p}\Biggl\{\frac{2}{n}\Biggl| \sum
_{i=1}^n\widehat{\xi}_{ir}\widetilde{\varepsilon}_i\Biggr|, \frac
{2}{n}\Biggl|\sum_{i=1}^n \widetilde{X}_{ij}\widetilde{\varepsilon}_i
\Biggr|\Biggr\} \le\frac{M_5\bolds{\alpha}_{\mathrm{sum}}}{v(k)^{3/2}}\sqrt
{\frac{\log p}{n}}.
\end{equation}
We will now prove (\ref{thmapp2-eq3}). For all $r=1,\ldots,k$ we have
\begin{eqnarray*}
\Biggl|\frac{1}{n}\sum_{i=1}^n\widehat{\xi}_{ir}\widetilde{\varepsilon
}_i\Biggr| &=& \Biggl|\frac{1}{n}\sum_{i=1}^n\widehat{\xi}_{ir}\sum
_{s=1}^k\alpha_s(\widehat{\xi}_{is}-\xi_{is})\Biggr|\\
&=& \Biggl|\frac{1}{n}\sum_{i=1}^n\widehat{\xi}_{ir}\sum_{s=1}^k\alpha
_s\biggl( \frac{\bolds{\psi}_s^T\mathbf{W}_i}{\sqrt{p\lambda_s}}
-\frac{\widehat{\bolds{\psi}}_s^T\mathbf{X}_i}{\sqrt{p\widehat
{\lambda}_s}} \biggr)\Biggr|\\
&=& \Biggl|\sum_{s=1}^k\alpha_s\frac{1}{n}\sum_{i=1}^n \widehat{\xi
}_{ir}\biggl(\frac{\bolds{\psi}_s^T\mathbf{X}_i}{\sqrt{p\lambda_s}}
- \frac{\widehat{\bolds{\psi}}_s^T\mathbf{X}_i}{\sqrt{p\widehat
{\lambda}_s}} - \frac{\bolds{\psi}_s^T\mathbf{Z}_i}{\sqrt{p\lambda
_s}}\biggr)\Biggr|\\
&\le& \alpha_{\mathrm{sum}} \sup_s\Biggl( \Biggl| \frac{1}{n}\sum_{i=1}^n \widehat
{\xi}_{ir}\frac{(\sqrt{\widehat{\lambda}_s}-\sqrt{\lambda_s})\widehat
{\bolds{\psi}}_s^T\mathbf{X}_i}{\sqrt{p\lambda_s\widehat{\lambda
}_s}}\Biggr| \\
&&\hspace*{44.3pt}{} + \Biggl|\frac{1}{n}\sum_{i=1}^n\widehat{\xi
}_{ir}\frac{(\bolds{\psi}_s-\widehat{\bolds{\psi}}_s)^T\mathbf
{X}_i}{\sqrt{p\lambda_s}}\Biggr| + \Biggl|\frac{1}{n}\sum_{i=1}^n\widehat
{\xi}_{ir}\frac{\bolds{\psi}_s^T\mathbf{Z}_i}{\sqrt{p\lambda
_s}}\Biggr|\Biggr).
\end{eqnarray*}
Using the Cauchy--Schwarz inequality and the fact that $\frac{1}{n}\sum
_{i=1}^n\widehat{\xi}_{ir}^2=1$, inequalities
(\ref{lem2-30}) and (\ref{lem2-31}) imply that under events
(\ref{B21})--(\ref{B24}),
one obtains
%
%
\begin{eqnarray} \label{thmapp2-eq4}\qquad
&&
\Biggl|\frac{1}{n}\sum_{i=1}^n\widehat{\xi}_{ir}\widetilde{\varepsilon
}_i\Biggr| \nonumber\\
&&\qquad\le \alpha_{\mathrm{sum}} \sup_s\Biggl(\Biggl( \frac{1}{n}\sum
_{i=1}^n\frac{(\sqrt{\widehat{\lambda}_s}-
\sqrt{\lambda_s})^2(\widehat{\bolds{\psi}}_s^T\mathbf
{X}_i)^2}{p\lambda_s\widehat{\lambda}_s}\Biggr)^{1/2}
\nonumber\\[-8pt]\\[-8pt]
&&\qquad\quad\hspace*{43.4pt}{} + \Biggl(\frac{1}{n}\sum_{i=1}^n\frac{((\bolds{\psi
}_s-\widehat{\bolds{\psi}}_s)^T\mathbf{X}_i)^2}{p\lambda_s}
\Biggr)^{1/2}+ \Biggl|\frac{1}{n}\sum_{i=1}^n\widehat{\xi}_{ir}\frac
{\bolds{\psi}_s^T\mathbf{Z}_i}{\sqrt{p\lambda_s}}\Biggr|
\Biggr)\nonumber\\
&&\qquad\leq \alpha_{\mathrm{sum}} \Biggl(\frac{\sqrt{M_3^{*}}}{v(k)}\sqrt{\frac{\log p}{n}}+
\frac{\sqrt{M_3^{**}}}{v(k)^{3/2}}\sqrt{\frac{\log p}{n}}
+\sup_s\Biggl|\frac{1}{n}\sum_{i=1}^n\widehat{\xi}_{ir}\frac{\bolds
{\psi}_s^T\mathbf{Z}_i}{\sqrt{p\lambda_s}}\Biggr|
\Biggr).\nonumber
\end{eqnarray}
Since also $\frac{1}{n}\sum_{i=1}^n\widetilde{X}_{ij}^2=1$, similar
arguments show that under (\ref{B21})--(\ref{B24})
%
%
\begin{eqnarray}
\Biggl|\frac{1}{n}\sum_{i=1}^n\widetilde{X}_{ij}\widetilde{\varepsilon
}_i\Biggr|
&\leq&\alpha_{\mathrm{sum}} \Biggl(\frac{\sqrt{M_3^{*}}}{v(k)}\sqrt{\frac{\log p}{n}}+
\frac{\sqrt{M_3^{**}}}{v(k)^{3/2}}\sqrt{\frac{\log
p}{n}}\nonumber\\[-8pt]\\[-8pt]
&&\hspace*{66.6pt}
{}+\sup_s\Biggl|\frac{1}{n}\sum_{i=1}^n\widetilde{X}_{ij}\frac{\bolds
{\psi}_s^T\mathbf{Z}_i}{\sqrt{p\lambda_s}}\Biggr|
\Biggr)\nonumber
\end{eqnarray}
for
all $j=1,\ldots,p$.
The Cauchy--Schwarz inequality yields $\sum_{l=1}^p |\widehat{\psi
}_{rl}|\leq\sqrt{p}$, $\sum_{l=1}^p |\psi_{rl}|\leq\sqrt{p}$, as well as
\begin{eqnarray*}
\sum_{l=1}^p \Biggl|\sum_{r=1}^k\widehat{\psi}_{rj}\widehat{\psi}_{rl}\Biggr|
&\le&\sqrt{p} \sqrt{\sum_{l=1}^p \Biggl(\sum_{r=1}^k\widehat{\psi}_{rj}\widehat
{\psi}_{rl}\Biggr)^2}
=\sqrt{p} \sqrt{ \sum_{r=1}^k\sum_{s=1}^k\widehat{\psi}_{rj}\widehat
{\psi}_{sj}\sum_{l=1}^p
\widehat{\psi}_{rl}\widehat{\psi}_{sl}}\\
&\le&{\sqrt{kp}\sup_r} |\widehat
{\psi}_{rj}|.
\end{eqnarray*}
Necessarily, $v(k)\le D_0/k$ and hence $k\le D_0/v(k)$.
It therefore follows from (\ref{thmapp1-3}), (\ref{thmapp1-4}),
(\ref{pkapp0}) and (A.5) that under events (\ref{B21})--(\ref{B24})
there are
some constants $M_5^{***}$, $\widetilde{M}_5^{***}$ such that for all
$r,s=1,\ldots,k$ and
$j=1,\ldots,p$,
%
%
\begin{eqnarray}\label{lem2-4}
\Biggl|\frac{1}{n}\sum_{i=1}^n\widehat{\xi}_{ir}\frac{\bolds{\psi
}_s^T\mathbf{Z}_i}
{\sqrt{p\lambda_s}}\Biggr|&=& \Biggl|\frac{1}{n}\sum_{i=1}^n\frac
{1}{p\sqrt{\widehat{\lambda}_r\lambda_s}}
\sum_{j=1}^p\sum_{j'=1}^p\widehat{\psi}_{rj}\psi
_{sj'}X_{ij}Z_{ij'}\Biggr|\nonumber\\[-8pt]\\[-8pt]
&\le&\frac{M_5^{***}}{v(k)}\sqrt{\frac{\log p}{n}}\nonumber
\end{eqnarray}
and
%
%
\begin{eqnarray}\label{lem2-5}\quad
\Biggl|\frac{1}{n}\sum_{i=1}^n\widetilde{X}_{ij}\frac{\bolds{\psi
}_s^T\mathbf{Z}_i}
{\sqrt{p\lambda_s}}\Biggr|&=&\Biggl|\frac{1}{n}\sum_{i=1}^n\sum
_{j'=1}^p\frac{(X_{ij}-
\sum_{r=1}^k\widehat{\psi}_{rj}\widehat{\bolds{\psi}}{}^T_r\mathbf
{X}_i)\psi_{sj'}Z_{ij'}}{( ({1}/{n})\sum_{i=1}^n(\widehat{\mathbf
{P}}_k\mathbf{X}_i)^2_j)^{1/2}\sqrt{p\lambda_s}}\Biggr|
\nonumber\\[-8pt]\\[-8pt]
&\le&\frac{\widetilde{M}_5^{***}}{v(k)^{3/2}}\sqrt{\frac{\log
p}{n}}.\nonumber
\end{eqnarray}
Result (\ref{thmapp2-eq3}) is now a direct consequence of (\ref
{thmapp2-eq4})--(\ref{lem2-5}). Note that all constants in
(\ref{thmapp2-eq4})--(\ref{lem2-5}) and thus also the constant
$M_5<\infty$ can be chosen independently of the values $n,p,k,S$
satisfying (A.3)--(A.5).

Under event (\ref{thmapp2-eq1}) as well as $\mathbf
{K}_{n,p}(k,S,3)>0$, inequalities (B.1), (4.1), (B.27) and (B.30) of
\citet{BiRiTss09} may be transferred in our context which yields
%
%
\begin{equation}
\|(\widehat{\bolds{\theta}}-\bolds{\theta})_{J_0}\|_2\le4\rho
\sqrt{k+S}/\kappa^2,\qquad \|\widehat{\bolds{\theta}}-\bolds{\theta
}\|_1\le4\|(\widehat{\bolds{\theta}}-\bolds{\theta})_{J_0}\|_1,
\end{equation}
where $J_0$ is the set of nonnull coefficients of $\bolds{\theta
}$. This implies that
%
%
\begin{equation}\label{bickel1}
\sum_{r=1}^k|\widehat{\widetilde{\alpha}}_r-\widetilde{\alpha}_r|+\sum
_{j=1}^p|\hspace*{1.5pt}\widehat{\hspace*{-1.5pt}\widetilde{\beta}}_j-\widetilde{\beta}_j|\le16\frac
{k+S}{\kappa^2}\rho.
\end{equation}

Events (\ref{B21})--(\ref{B24}) hold with probability $A(n,p)$, and
therefore the probability of event (\ref{thmapp2-eq1}) is at least
$A(n,p)-(p+k)^{1-A^2/8}$. When combining
(\ref{thmapp1-4}), (\ref{pkapp0}) and (\ref{bickel1}), inequalities
(\ref{thmaug-1}) and (\ref{thmaug-2}) follow from
the definitions of $\widetilde{\beta}_j$ and $\widetilde{\alpha}_r$,
since under
(\ref{B21})--(\ref{B24})
\begin{eqnarray*}
\sum_{j=1}^p| \widehat{\beta}_j-\beta_j | &=& \sum_{j=1}^p\frac
{|\hspace*{1.5pt}\widehat{\hspace*{-1.5pt}\widetilde{\beta}}_j
-\widetilde{\beta}_j|}{(
({1}/{n})\sum_{i=1}^n (\widehat{\mathbf{P}}_k\mathbf{X}_i)_j^2)^{1/2}}
\\
&\le&\frac{\sum_{j=1}^p|\hspace*{1.5pt}\widehat{\hspace*{-1.5pt}\widetilde{\beta}}_j-\widetilde{\beta
}_j|}{(D_1-M_1({kn^{-1/2}\sqrt{\log p}}/{v(k)^{1/2}}))^{1/2}}
\end{eqnarray*}
and
\begin{eqnarray*}
\sum_{r=1}^k|\widehat{\alpha}_r-\alpha_r| &=& \sum_{r=1}^k\Biggl|\widehat
{\widetilde{\alpha}}_r-\widetilde{\alpha}_r-\sqrt{p\widehat{\lambda
}_r}\sum_{j=1}^p\widehat{\psi}_{rj}(\widehat{\beta}_j-\beta_j)\Biggr|\\
&\le& \sum_{r=1}^k|\widehat{\widetilde{\alpha}}_r-\widetilde{\alpha
}_r|+ k\bigl(D_0+C_0n^{-1/2}\sqrt{\log p}\bigr)^{1/2}\sum_{j=1}^p|\widehat{\beta
}_j-\beta_j|.
\end{eqnarray*}
It remains to prove assertion (\ref{thmaug2-1}) on the prediction
error. We have
%
%
\begin{eqnarray}\label{thmaug2-eq1}\qquad
&&\frac{1}{n}\sum_{i=1}^n\Biggl(\sum_{r=1}^k\widehat{\xi}_{ir}\widehat
{\alpha}_r-\xi_{ir}\alpha_r+\sum_{j=1}^pX_{ij}(\widehat{\beta}_j-\beta
_j)\Biggr)^2\nonumber\\
&&\qquad=\frac{1}{n}\sum_{i=1}^n\Biggl(\sum_{r=1}^k\widehat{\xi}_{ir}(\widehat
{\widetilde{\alpha}}_r-\widetilde{\alpha}_r)+\sum_{j=1}^p\widetilde{X}_{ij}
(\hspace*{1.5pt}\widehat{\hspace*{-1.5pt}\widetilde{\beta}}_j-\widetilde{\beta}_j)+\sum
_{r=1}^k(\widehat{\xi}_{ir}-\xi_{ir})\alpha_r\Biggr)^2\nonumber\\[-8pt]\\[-8pt]
&&\qquad\le\frac{2}{n}\sum_{i=1}^n\Biggl(\sum_{r=1}^k\widehat{\xi
}_{ir}(\widehat{\widetilde{\alpha}}_r-\widetilde{\alpha}_r)
+\sum_{j=1}^p\widetilde{X}_{ij}(\hspace*{1.5pt}\widehat{\hspace*{-1.5pt}\widetilde{\beta}}_j-\widetilde
{\beta}_j)\Biggr)^2\nonumber\\
&&\qquad\quad{}+
\frac{2}{n}\sum_{i=1}^n\Biggl(\sum_{r=1}^k(\widehat{\xi}_{ir}-\xi
_{ir})\alpha_r\Biggr)^2.\nonumber
\end{eqnarray}
Under event (\ref{thmapp2-eq1}) as well as $\mathbf
{K}_{n,p}(k,S,3)>0$, the first part of inequalities (B.31) in the proof of
Theorem 7.2 of \citet{BiRiTss09} leads to
%
%
\begin{equation}\label{thmaug2-eq2}
\frac{2}{n}\sum_{i=1}^n\Biggl(\sum_{r=1}^k\widehat{\xi}_{ir}(\widehat
{\widetilde{\alpha}}_r-\widetilde{\alpha}_r)+\sum_{j=1}^p\widetilde
{X}_{ij}(\hspace*{1.5pt}\widehat{\hspace*{-1.5pt}\widetilde{\beta}}_j-\widetilde{\beta}_j)\Biggr)^2\le
\frac{32(k+S)}{\kappa^2}\rho^2.
\end{equation}
Under events (\ref{B21})--(\ref{B24}), (\ref{thmapp2-eq1}) as well as
(\ref{B25}), inequality (\ref{thmaug2-1}) now follows from (\ref
{thmaug2-eq1}), (\ref{thmaug2-eq2}) and (\ref{pkapp2}). The assertion
then is a consequence
of the fact that (\ref{B21})--(\ref{B24}) are satisfied with
probability $A(n,p)$, while (\ref{thmapp2-eq1}) and
(\ref{B25}) hold with probabilities at least $1-(p+k)^{-A^2/8}$ and
$1-\frac{k}{n}$, respectively.
\end{pf*}
\end{appendix}


%

\printaddresses

\end{document}